\documentclass{amsart}
\usepackage{amscd}
\usepackage{amssymb}
\input xypic
\xyoption{all}
\newdir{ >}{{}*!/-5pt/@{>}}

\newcommand{\en}{\operatorname{End}}

\renewcommand{\r}{\rightarrow}
\newcommand{\kh}{\operatorname{KH}}
\newcommand{\alg}{\operatorname{Alg}}

\newcommand{\coli}{\operatorname{colim}}
\newcommand{\limi}{\operatorname{lim\mbox{ }}}

\def\id{\operatorname{id}}

\def\cC{\mathcal C}
\def\cV{\mathcal V}
\def\cT{\mathcal T}

\def\cD{\mathcal D}

\def\cE{\mathcal E}

\def\cU{\mathcal U}

\def\cF{\mathcal F}

\def\cL{\mathcal L}
\def\cM{\mathcal M}

\def\cP{\mathcal P}
\def\cS{\mathcal S}
\def\cK{\mathcal K}

\def\ev{\mathrm{ev}}

\newcommand{\C}{\mathbb{C}}

\newcommand{\Z}{\mathbb{Z}}
\newcommand{\N}{\mathbb{N}}
\newcommand{\kk}{{\operatorname{kk}}}
\newcommand{\sd}{{\operatorname{sd}}}
\newcommand{\hkk}{{\operatorname{kk}}^{G}}
\newcommand{\indg}{\operatorname{Ind^{G}_{H}}}

\newcommand{\resg}{\operatorname{Res_{G}^{H}}}

\newcommand{\sS}{\mathbb{S}}

\newcommand{\Cat}{\mathcal{D}}

\def\map{\operatorname{map}}

\numberwithin{equation}{section}
\theoremstyle{plain}
\newtheorem{thm}[equation]{Theorem}
\newtheorem{cor}[equation]{Corollary}

\newtheorem{prop}[equation]{Proposition}

\newcommand{\comment}[1]{}  

\theoremstyle{definition}

\theoremstyle{remark}
\newtheorem{rem}[equation]{Remark}
\newtheorem{ex}[equation]{Example}

\begin{document}

\bibliographystyle{plain}

\title{Equivariant algebraic kk-theory and adjointness theorems}
\author{Eugenia Ellis}
\email{eugenia@cmat.edu.uy - eugenia.ellis@fu-berlin.de}
\address{CMAT, Facultad de Ciencias-UDELAR\\
Igu\'a $4225$, 11400 Montevideo, Uruguay \newline Institut f\"ur Mathematik, FU, Arnimallee 7, Raum 212, 14195 Berlin, Germany.}
\thanks{The author was partially supported by CSIC-Uruguay; by the MathAmSud network U11MATH-05 partially financed by ANII (Uruguay) and MINCyT (Argentina) and by grants UBACyT 20020100100386, and MTM2007-64704 (FEDER funds).}

\begin{abstract}
We introduce an equivariant algebraic $kk$-theory for $G$-algebras and $G$-graded algebras. We study some adjointness theorems related with crossed product, trivial action, induction and restriction. In particular we obtain an algebraic version of the Green-Julg Theorem which gives us a computational tool.  
\end{abstract}

\maketitle
\section{Introduction}

Algebraic $kk$-theory has been introduced by  G. Corti\~nas  and A. Thom
in \cite{CT}. This is a bivariant $K$-theory on the category of
$\ell$-algebras where $\ell$ is a commutative ring with unit.
For each pair $(A,B)$ of $\ell$-algebras a group $\kk(A,B)$ is defined. 
A category $\mathfrak{KK}$ is obtained whose objects are 
$\ell$-algebras and where the morphisms from $A$ to $B$ are the elements of the group $\kk(A,B)$. The category $\mathfrak{KK}$ is triangulated and there is
a canonical functor $j: \alg_{\ell}\r \mathfrak{KK}$ with universal
properties. These properties are algebraic homotopy invariance, matrix
invariance and excision.

The definition of algebraic $kk$-theory was inspired by the work
of J. Cuntz \cite{ckkalg} and N. Higson \cite{hkk} on the universal
properties of Kasparov 
$KK$-theory \cite{kkk}. The $KK$-theory of separable $C^{*}$-algebras is a
common generalization both of topological $K$-homology and topological
$K$-theory as an
additive bivariant functor. Let $A$, $B$ be separable $C^{*}$-algebras. Then
\begin{equation}\label{kkgen2} 
KK_{*}(\C,B)\simeq K^{top}_{*}(B)  \qquad KK^{*}(A,\C)= {K^{*}_{hom}(A)} 
\end{equation}       
here $K^{top}_{*}(B)$ denotes the $K$-theory of $B$ and
$K^{{*}}_{hom}(A)$ the topological K-homology of $A$.
J. Cuntz in \cite{ckk} gave another equivalent definition of the
original one given in \cite{kkk}. This new approach allowed to put
bivariant $K$-theory in
algebraic context. Higson in \cite{hkk}
stated the universal property of $KK$ whose algebraic analogue
is studied in \cite{CT}, where also an analogue of
\eqref{kkgen2} is proved. On the algebraic side, if $A$ is an
$\ell$-algebra then  
$$
\kk(\ell, A) \simeq \kh (A)
$$  
here $\kh$ is Weibel's homotopy K-theory defined in \cite{kh}.
We can start to build a dictionary between Kasparov's KK-theory and
algebraic kk-theory in the following way
{\footnotesize{
\mbox{}
\begin{center}
\begin{tabular}{|c c c|}
\hline
&&\\
Kasparov's $KK$-theory & $\leftrightarrow$  & algebraic $kk$-theory \\
&&\\
bivariant $K$-theory on & & bivariant $K$-theory on \\
separable $C^{*}$-algebras & $\leftrightarrow$& $\ell$-algebras \\
$C^{*}$-$\alg$&& $\alg_{\ell}$\\
&&\\
 $k: C^{*}$-$\alg \r {KK}$ & $\leftrightarrow$&
 $j: \alg_{\ell} \r \mathfrak{KK}$\\
&&\\
$k$ is stable with respect to  & $\leftrightarrow$ & $j$ is stable with respect to \\
compact operators $\cK$ && $M_{\infty}=\bigcup_{n\in \N} M_{n}$\\
&&\\
$k$ is continous homotopy invariant & $\leftrightarrow$ & $j$ is
polynomial homotopy invariant \\
&&\\
$k$ is split exact & $\leftrightarrow$& $j$ is excisive \\ 
&&\\
$k$ is universal for the properties & $\leftrightarrow$ & $j$ is
universal for the 
properties\\ 
described above && described above \\
&&\\
$\operatorname{KK}_{*}(\C,A)\simeq K^{top}_{*}(A)$ & $\leftrightarrow$ & $\kk_{*}(\ell, A) \simeq \kh_{*}(A)$\\
\hline
\end{tabular}
\end{center} 
}}
\vspace{0.2cm}
In this paper we obtain an equivariant version of the dictionary stated above in the following sense
\vspace{0.2cm}
\mbox{}
{\footnotesize
\begin{center}
\begin{tabular}{|c c c|}
\hline&&\\
Equivariant Kasparov's $KK$-theory & $\leftrightarrow$ & Equivariant algebraic $kk$-theory \\
&&\\
bivariant $K$-theory on && bivariant $K$-theory on \\
separable $G$-$C^{*}$-algebras & $\leftrightarrow$& $G$-algebras \\
&&\\
 $k:G$-$C^{*}$-$\alg \r {KK}^{G}$ &$\leftrightarrow$ &
 $j^{G}: G$-$\alg \r \mathfrak{KK}^{G}$\\
&&\\
$k$ is stable with respect to  &  $\leftrightarrow$& $j^{G}$ is $G$-stable \\
$\cK(\ell^{2}(G\times \N))$& &\\
&&\\
$k$ is continous homotopy invariant & $\leftrightarrow$& $j^{G}$ is polynomial homotopy invariant \\
&&\\
$k$ is split exact & $\leftrightarrow$& $j^{G}$ is excisive \\ 
&&\\
$k$ is universal for the properties & $\leftrightarrow$& $j^{G}$ is universal for the
properties\\
described above& & described above \\
&&\\
$\operatorname{KK^{G}}_{*}(\C,A)\simeq K^{top}_{*}(A\rtimes G)$ & $\leftrightarrow$ & $\kk^{G}(\ell, A) \simeq \kh (A\rtimes G)$\\
with $G$ compact & & with $G$ finite and $1/|G| \in \ell$\\
\hline
\end{tabular}
\end{center}
}
We also introduce a dual theory $\hat{\mathfrak{KK}}^{G}$ for $G$-graded algebras and establish a duality result similar to that proved by Baaj and Skandalis in 
\cite{bs}.

In Section \ref{akk} we recall some results from \cite{CT}. 
We take special care in the definition of equivariant matrix invariance.  We introduce the concept of $G$-stable functor in Section \ref{eqmi}, the rest of the concepts which appears in the equivariant dictionary are straightforward. The definition of $G$-stability was inspired by the definition of
 equivariant stability for $G$-$C^{*}$-algebras (see \cite{Mkk}). 
In sections \ref{akk} and \ref{eqmi} $G$ is a group without any other assumption but from Section \ref{ekk} and for the rest of the paper $G$ is a countable group. 

In Section \ref{ekk} we introduce the appropiate brand of equivariant algebraic
$kk$-theory in each case and establish its universal properties. 
For a countable group $G$, we define an equivariant algebraic $kk$-theory  $\mathfrak{KK}^{G}$ for the category of $G$-algebras and
 $\hat{\mathfrak{KK}}^{G}$ for the category of $G$-graded algebras.  

We study adjointness theorems in equivariant $kk$-theory. 
We put in an algebraic context some of the adjointness theorems which
appear in Kasparov KK-theory. 
In Section \ref{tricrosec}, we define the functors of trivial action
and crossed product between $\mathfrak{KK}$ and
$\mathfrak{KK}^{G}$. The first adjointness theorem is Theorem
\ref{gjgalg} which is an algebraic version of the {\sc Green-Julg Theorem}. 
This result gives us the first computation related with 
homotopy K-theory. If $G$ is a
finite group, $A$ is a $G$-algebra, $B$ is an algebra and $\frac{1}{|G|}\in \ell$ then
there is an isomorphism
$$
\psi_{GJ}:\kk^{G}(B^{\tau},A)\r \kk(B, A \rtimes G).
$$
In particular, if $B=\ell$ then
$$
\kk^{G}(\ell,A)\simeq \kh(A\rtimes G).
$$ 

In Section \ref{inressec} we consider $H$ a subgroup of $G$. We define
induction and restriction funtors between $\mathfrak{KK}^{G}$ and
$\mathfrak{KK}^{H}$ and study the
adjointness between them. 
If $B$ is an $H$-algebra and $A$ is a $G$-algebra then there is an isomorphism 
$$
\psi_{IR}: \hkk(\indg B,A) \r {\kk}^{H}(B,\resg A).
$$
This result gives us another computation. Taking $H$ the trivial group and $B=\ell$ we obtain that 
$$
 \kk^{G}(\ell^{(G)},A)\simeq \kh(A) \qquad \forall A\in G-\alg.
$$
Here 
$\ell^{(G)}=\bigoplus_{g\in G} \ell$ with the regular action of $G$. 
More general, if $H$ is a finite subgroup of $G$ and $1/|H| \in \ell$ we combine $\psi_{GJ}$ and $\psi_{IR}$ and obtain  
$$
 \kk^{G}(\ell^{(G/H)},A)\simeq \kh(A\rtimes H) \qquad \forall A\in G-\alg.
$$  

In Section \ref{dualsec} we obtain an algebraic version of the Baaj-Skandalis duality theorem. We
show that the functors 
$$
\rtimes G: \mathfrak{KK}^{G}\rightarrow \hat{\mathfrak{KK}}^{G} \qquad \quad 
G \hat{\rtimes}:\hat{\mathfrak{KK}}^{G}\rightarrow \mathfrak{KK}^{G}
$$
are inverse category equivalences.
\subsection*{Acknowledgements} The results here are part of my PhD Thesis. 
I am thankful to Willie Corti\~nas for his orientation and comments about this paper. I would like also to thank the referee for useful comments. 
This paper was finished during a visit to the Freie
Universit\"at Berlin (FU), financed by Berlin Mathematical School (BMS) and International Mathematical Union (IMU). 
I am also thankful to Freie Universit\"at for its hospitality and to BMS-IMU
for its support.

\section{Algebraic kk-theory}\label{akk}
\numberwithin{equation}{subsection}
In this section we recall some results from \cite{CT} and we adapt them to
our setting.
Let $\ell$ be a commutative ring with unit. We consider an $\ell$-bimodule $A$ such that $x\cdot a = a\cdot x$, with $x\in \ell$ and $a\in A$.
An $\ell$-bimodule $A$ is an {\sf{$\ell$-algebra}} if it is an associative and not necessarily unital algebra.
Let $G$ be a group. A {\sf $G$-algebra} $A$ is an $\ell$-algebra with an action of $G$, i.e. with a group homomorphism $\alpha: G\rightarrow \en_{\ell}(A)$, where $\en_{\ell}(A)$ denote the group of $\ell$-linear endomorphisms of $A$. We shall denote by $g(a)$ or $g\cdot a$ the element $\alpha(g)(a)$. An {\sf equivariant morphism} $f:A\rightarrow B$ between $G$-algebras is a $G$-equivariant  $\ell$-linear map.
A {\sf{$G$-graded algebra}} $A$ is an $\ell$-algebra with a
familiy of $\ell$-submodules $\{A_{s}\}_{s\in G}$ such that
$$
A=\bigoplus_{s\in G}A_{s}\qquad\qquad 
A_{s}A_{t}\subseteq A_{st} \qquad\qquad s,t\in G. 
$$
We write $|a|=s$ if $a\in A_{s}$.
An {\sf homogeneous morphism} $f:A\r B$ of $G$-graded algebras is an algebra morphism  such that $f(A_{s})\subset B_{s}$, for all $s\in G$.
We consider the category of $G$-algebras with equivariant morphisms and the category of $G$-graded algebras with homogeneous morphisms. 

In this section we define an algebraic $kk$-theory for the categories of $G$-algebras and
$G$-graded algebras. We write $\cC$ to refer to either of these categories.  

\subsection{Homotopy invariance}
Let $A$ be an object of $\cC$. Put $A^{\Delta^{1}}:=A[t]=A\otimes_{\Z}
\Z[t].$ Consider the trivial action of $G$ (trivial $G$-grading) on $\Z[t]$ and
the diagonal structure on $A[t]$. Then $A[t]$ is an object of $\cC$.
Let us write $c_{A}:A\r A[t]$ for the inclusion of $A$ as constant polynomials  in $A[t]$ and $\ev_{i}:A[t]\r A$ for the evaluation of $t$ at $i$ 
($i=0,1$). Note that these maps are morphisms in $\cC$ and that $c_{A}$
is a section of $\ev_{i}$.

Let $f_{0},f_{1}:A\r B$ be morphisms in $\cC$. We call 
$f_{0}$ and $f_{1}$ {\sf elementarily homotopic}, and write $f_{0}\sim_{e} f_{1}$, if there exists a morphism $H:A\r B[t]$ such that 
 $\ev_{i}H=f_{i}$, $i=0,1$. 
It is easy to check that elementary homotopy is a reflexive and
symmetric relation but in general it is not transitive. 
Let $f,g:A\r B$ be morphisms in $\cC$. 
We call $f$ and $g$ {\sf homotopic}, and write $f \sim g$, if
they can be connected by a finite chain of elementary homotopies,
$$
f\sim_{e}h_{0}\sim_{e}\ldots\sim_{e}h_{n}\sim_{e}g.
$$
We denote the set of homotopy clases by $[A,B]_{\cC}$.
A morphism $f:A\r B$ is an {\sf elementary homotopy equivalence}
if there exists a morphism $g:B\r A$ such that $f\circ g\sim_{e} \operatorname{id}_{B}$ and $g\circ f\sim_{e}\operatorname{id}_{A}$. 
We say $A$ is {\sf elementarily  
contractible} if the null morphism and the identity morphism are elementarily homotopic.

The {\sf category of $\operatorname{ind}$-objects} of 
$\cC$ is the category $\operatorname{ind}$-$\cC$ of directed diagrams in $\cC$.
An  object in $\operatorname{ind}$-$\cC$ is described by a  
filtering partially ordered set $(I,\leq)$ and a functor $A:I \r \cC$.
The set of homomorphisms in $\operatorname{ind}$-$\cC$ is defined by 
$$
\displaystyle{
\hom_{\mbox{$\operatorname{ind}$-$\cC$}}((A,I),(B,J)):= \underset{i\in I}{\limi}\underset{j\in J}{\coli}
\hom_{\cC}(A_{i},B_{j}).
}
$$
Let $A=(A,I)$ and $B=(B,J)$ be objects of $\operatorname{ind}$-$\cC$, we have a map
\begin{equation}\label{homin}
\hom_{\mbox{$\operatorname{ind}$-$\cC$}}(A,B)\mapsto [A,B]_{\cC}=\underset{i\in I}{\limi}\underset{j\in J}{\coli}[A_{i},B_{j}]_{\cC}.
\end{equation}
We say two morphisms in $\operatorname{ind}$-${\cC}$ are {\sf
homotopic} if their images by \eqref{homin} are equal. An object of
$\operatorname{ind}$-$\cC$ is {\sf{contractible}} if the null morphism
and the identity morphism are homotopic.
Sometimes we shall omit the poset in the notation. 

Let $\Cat$ be an arbitrary category.
A functor $F:\cC \r \Cat$ is  {\sf homotopy invariant} if it maps the inclusion $c_{A}: A\r A[t]$ to an isomorphism. 
It is easy to check that $F:\cC \r \Cat$ is an homotopy invariant functor if and only if $F(f)=F(g)$ when $f\sim g$ (or equivalenty when $f\sim_{e} g$).  If $\cC$ is the category of $G$-algebras we say that an homotopy invariant functor is {\sf equivariantly homotopy invariant} and if $\cC$ is the category of $G$-graded algebras we say that a homotopy invariant functor is {\sf graded homotopy invariant}.

\subsection{Matrix stability}\label{maes}
Consider $M_{n}$ the algebra of $n\times n$-matrices with coefficients in $\Z$ with the trivial action (grading) of $G$ and $M_{\infty}=\cup_{n\geq 0} M_{n}$.
Let $A$ be an object in $\cC$. We define 
$$
M_{n} A = M_{n} \otimes_{\Z} A \qquad  M_{\infty} A =
M_{\infty} \otimes_{\Z} A
$$ 
which are objects in $\cC$ with the diagonal action (grading). 
Denote by $\iota_{n}: \Z \r M_{n}$ and $\iota_{\infty}: \Z \r
M_{\infty}$ the inclusions at the upper left corner.
A functor $F:\cC \r \Cat$ is {\sf $M_{n}$-stable} 
{(\sf $M_{\infty}$-stable)} if $F(\iota_{n}\otimes \id_{A})$
($F(\iota_{\infty}\otimes \id_{A})$) is an isomorphism for all 
$A \in \cC$.

\subsection{Algebra of polynomial functions}
Consider the following simplicial ring
\begin{equation}\label{defalsi}
\Z^{\Delta}: [n]\mapsto \Z^{\Delta^{n}} \quad\quad
\Z^{\Delta^{n}}:=\Z[t_{0},\ldots, t_{n}]/<1-\sum_{i}t_{i}>
\end{equation}
$$
\begin{array}{lll}
\Theta: [n]\r [m] & \mapsto & {\Theta^{*}}: \Z^{\Delta^{m}}\r
\Z^{\Delta^{n}}\\
& &{\Theta^{*}}( t_{i}) =\left\{ \begin{array}{ll}0 & \mbox{ si }
\Theta^{-1}(i)=\emptyset
\\
\sum_{j\in \Theta^{-1}(i)}t_{j} & \mbox{ si } \Theta^{-1}(i)\neq \emptyset
 \end{array}\right. 
\end{array}
$$
Let $A$ be an object in $\cC$. Define
$$
A^{\Delta}: [n]\mapsto A^{\Delta^{n}} \quad\quad
A^{\Delta^{n}}:= A\otimes_{\Z}\Z^{\Delta^{n}}
$$
Note $A^{\Delta}$ is a simplicial $\cC$-object. 
Let  $X$ be a simplicial set. Define 
$$
A^{X}:= \map_{\sS}(X,A^{\Delta})
$$
If $(K,\star)$ is a pointed simplicial set, put
$$
A^{(K,\star)}:=\map_{\sS_{*}}((K,\star),A^{\Delta})=
\ker(\map_{\sS}(K,A^{\Delta})\r \map_{\sS}(\star,A^{\Delta}))=\ker(A^{K}\r A)
$$
If $A$ is a $G$-algebra ($G$-graded algebra) and $K$ a finite
simplicial set, by \cite[Lemma 3.1.3]{CT} we can consider $A^{K}$ and $A^{(K,\star)}$ as  a
$G$-algebras ($G$-graded algebras) taking the diagonal structure with the trivial structure in $\Z^{K}$ and $\Z^{(K,\star)}$.

We will denote by $\sd^{\bullet} X$ the following pro-simplicial set
$$
\sd^{\bullet} X: \ldots \xrightarrow{} \sd^{n}X\xrightarrow{h} \sd^{n-1}X
\xrightarrow{}\ldots \xrightarrow{} \sd X \xrightarrow{h} X
$$ 
where $\sd X$ is the subdivision of $X$ and $h$ is the last vertex map, see \cite[III.4]{GoJ}.  
If $A$ is an object in $\cC$ we consider the following  object in $\operatorname{ind}$-$\cC$
$$
A^{\sd^{\bullet} X}: A^{X}\r A^{\sd X} \r\ldots\r
A^{\sd^{n-1}X}\r  A^{\sd^{n}X} \r \ldots 
$$ 
\subsection{Extensions and classifying map}
A sequence of morphisms in $\operatorname{ind}$-$\cC$ 
\begin{equation}\label{exten}
A\xrightarrow{f} B \xrightarrow{g} C
\end{equation}
is called an  {\sf extension} if $f$ is a kernel of $g$ and $g$ is a
cokernel of $f$.
Let $\mathcal{U}(\cC)$ be the category of modules with linear and
equivariant maps (graded maps), i.e. we forget the multiplication in $\cC$ keeping the module structure. Let $F:\cC \r \mbox{$\mathcal{U}(\cC)$}$ be the forgetful functor. 
This functor can be extended to 
$F:\operatorname{ind}$-$\cC\r \operatorname{ind}$-$\mathcal{U}(\cC)$.
We will call an extension \eqref{exten}  {\sf
  weakly split}  if $F(g)$ has a section in $\mbox{ $\operatorname{ind}$-$\mathcal{U}(\cC)$}$. 

Let $M$ be an object in $\cU(\cC)$. Consider in
$$
\tilde{T}(M)=\bigoplus_{n\geq 1} M^{\otimes^{n}}\quad
M^{\otimes^{n}}=\underset{\mbox{$n$-times}}{\underbrace{M\otimes
\ldots \otimes M}}
$$
the usual structure of $\ell$-algebra. If $\cC$ is the category of $G$-algebras, we consider in $M^{\otimes^{n}}$ the following action,
$$
g\cdot(m_{1}\otimes m_{2} \otimes \ldots \otimes m_{n})= 
g\cdot m_{1}\otimes g\cdot m_{2} \otimes \ldots \otimes g\cdot m_{n}
$$
which gives to $\tilde{T}(M)$ a $G$-algebra structure.  If $\cC$ is the category of $G$-graded algebras, take in $M^{\otimes^{n}}$ the following $G$-graduation
$$
|m_{1}\otimes m_{2} \otimes \ldots \otimes m_{n}|=  |m_{1}||m_{2}| \ldots|m_{n}|.
$$
By this way $\tilde{T}(M)$ is a $G$-graded algebra. 
Both constructions are functorial hence we consider the functor
 $\tilde{T}:\cU(\cC)\r \cC$. 
Put
$$
T:=\tilde{T}\circ F:\cC \r \cC
$$ 
If $A$ is an object in $\cC$ there exists an morphism in $\cC$
$$
\eta_{A}: T(A)\r A \qquad \eta_{A}(a_{1}\otimes \ldots \otimes
a_{n})=a_{1}\ldots a_{n}
$$
and a morphism in $\cU(\cC)$ $\mu_{A}:A\r T(A)$ which is the inclusion at the first summand of $T(A)$.
Let $A,B$ be objects in $\cC$, it is easy to check that
\begin{equation}\label{adjadj}
\hom_{\cC}(T(A), B) \simeq
\hom_{\cU(\cC)}(F(A),F(B)).
\end{equation}
Hence if we have a morphism $A\r B$ in $\cU(\cC)$, we can
extend it to a morphism $T(A)\r B$ in $\cC$. 
It shows that $T$ is the left adjoint of $F$.  
The counit of the adjuntion is $\eta_{A}: T(A)\r A$ and it is surjective (see
\cite[IV.3 Thm 1]{Mac}). We define 
$$ J(A):=\ker \eta_{A}.$$
The {\sf universal extension} of $A$ is
$$
J(A)\xrightarrow{\imath_{A}} T(A)\xrightarrow{\eta_{A}} A.
$$
Let $A\xrightarrow{f} B \xrightarrow{g} C$ be a weakly split extension. Let $s$ be a section of $F(g)$ and define $\hat{\xi}=\eta_{B}\circ
\tilde{T}(s)$, then
$$
\eta_{C}=\eta_{C}\circ T(g)\circ \tilde{T}(s)= g \circ \eta_{B}\circ
\tilde{T}(s)= g\circ \hat{\xi}.
$$
Define $\xi: J(C)\r A$ as the restriction of $\hat{\xi}$ to $J(C)$
We obtain a commutative diagram of extensions
$$
\xymatrix{
A\ar[r]^-{f}& B\ar[r]^-{g}& C\\
J(C)\ar[u]^-{\xi}\ar[r]_-{\imath_{C}}& T(C)\ar[u]^-{\hat{\xi}}\ar[r]_-{\eta_{C}}&C\ar[u]_-{\id_{C}}
}
$$
This
construction of $\xi$ depends on which $s$ is chosen.
If $\xi_{1},\xi_{2}:J(C)\r A$ are morphisms constructed taking
different sections of  $F(g)$ we will show $\xi_{1},\xi_{2}$ are homotopic.
Define a linear equivariant (graded) map
$$
H: C\r A[t] \quad H(c)= (1-t)\xi_{1}(c) + t \xi_{2}(c).
$$
By \eqref{adjadj} it extends to a homomorphism and there exists a morphism $H:T(C)\r A[t]$ in $\cC$ such that 
$$
\ev_{0}\circ H_{|J(C)}= \xi_{1} \qquad \qquad \ev_{1}\circ H_{|J(C)}= \xi_{2}.
$$
Then the map $\xi$ is unique up to elementary homotopy. We call $\xi$ {\sf the classifying map} of the extension $A\xrightarrow{f} B \xrightarrow{g} C$.
This construction is functorial, that means if we have a diagram of  weakly split extensions,
$$
\xymatrix{
A\ar[r]\ar[d]_{f} & B \ar[d]^{h}\ar[r]& C\ar[d]^{g}\\
A'\ar[r] &B'\ar[r]& C'
}
$$
then there is a
diagram 
$$
\xymatrix{
J(C)\ar[r]\ar[d]_{J(g)}& A\ar[d]^{f}\\
J(C')\ar[r]& A'
}
$$
of classifying maps, which is commutative up to elementary homotopy.

Let $L$ be a ring and $A$ an object in $\cC$. Then the extension  
\begin{equation}\label{jaltal}
J(A)\otimes_{\Z} L\r T(A)\otimes_{\Z}L\r A\otimes_{\Z}L
\end{equation}
is weakly split, and there is a choice for classifying map $$\phi_{A,L}:
J(A\otimes_{\Z}L)\r J(A)\otimes_{\Z}L$$ of \eqref{jaltal}, which is
natural in both variables. In particular, if $K$ is a finite pointed simplicial set. There is a homotopy class
of maps 
\begin{equation}\label{jaka}
\phi_{A,K}:J(A^{K})\r J(A)^{K}
\end{equation}
natural with respect to $K$, which is represented by a classifying map of the following extension 
$$
J(A)^{K}\xrightarrow{\iota^{*}_{A}} T(A)^{K}\xrightarrow{\eta_{A}^{*}} A^{K}
$$

Let $S^{1}$ be the simplicial circle $\Delta^{1}/\partial \Delta^{1}$,
we define 
$$
\Omega :=\Z^{(S^{1},\star)} \qquad P:= \Z^{(\Delta^{1},\star)}
$$
The {\sf path extension} of $A$ is the extension induced by the cofibration 
$\partial \Delta^{1}\subset \Delta^{1}$ (see \cite[Lemma 3.1.2]{CT}),
\begin{equation}\label{excam}
\Omega A \xrightarrow{\quad\qquad} P A \xrightarrow{(\ev_{0},\ev_{1})} A\oplus A 
\end{equation}
The extension \eqref{excam} is weakly split because we have a
 linear and equivariant (graded) section of $(\ev_{0},\ev_{1})$
\begin{equation}\label{secec}
(a,b)\mapsto (1-t)a + tb.
\end{equation}
The {\sf loop extension} of $A$ is 
\begin{equation}{\label{extlaz}}
\Omega A\r PA\xrightarrow{\ev_{1}}A
\end{equation}
Note $a\mapsto at$ is a natural section of $F(\ev_{1})$. Thus we can pick a natural choice for the classifying map of
\eqref{extlaz}. We call it
$$
\rho_{A}:J(A)\r \Omega A
$$

Let $\cP:=\Z^{(\sd^{\bullet}\Delta^{1},\star)}$
we have an extension in $\operatorname{ind}$-$\cC$
$$
A^{\cS^{1}}:=A^{(\sd^{\bullet} S^{1},\star)} \quad\quad
A^{\cS^{1}} \r \cP
A\xrightarrow{\ev_{1}} A
$$
which is naturally weakly split. The classifying map $J(A)\r
A^{\cS^{1}}$ is the following composition
\begin{equation}{\label{morcla}}
J(A)\xrightarrow{\rho_{A}}\Omega A\xrightarrow{h}A^{\cS^{1}} 
\end{equation}
where $h$ is induced by the last vertex map. We will sometimes abuse
notation and write $\rho_{A}$ for the map \eqref{morcla}.

Let $f:A\r B$ be a morphism in $\cC$. The {\sf mapping
path extension of $f$} is the extension obtained from the loop extension of
$B$ by pulling it back to $A$
$$
P_{f}:= PB \times_{B}A\quad\quad \xymatrix{
\Omega B \ar[r]^-{\iota}\ar[d] & PB \times_{B}
A\ar[r]^-{\pi_{f}}\ar[d] & A\ar[d]^{f}& \operatorname{E}'\\
\Omega B\ar[r]_-{\iota}& \ar[r] PB \ar[r]_-{\ev_{1}}& B& \operatorname{E}
}
$$
We call $P_{f}$ the {\sf path algebra} of $f$. 
Note the extension $\operatorname{E}'$ is naturally weakly split
because $\tilde{s}(a)=(s\circ f(a),a)$ is a natural section of
$\pi_{f}$ where $s$ is the natural section of $\operatorname{E}$.
We define
\begin{equation}\label{extenpf}
\cP_{f}:=\cP B\times_{B} A \quad \quad B^{\cS^{1}} \xrightarrow{} \cP
B\times_{B} A \xrightarrow{\pi_{f}} A.
\end{equation}

Note $\rho_{f}:=\rho_{B}J(f)$ is the classifying map of the extension \eqref{extenpf}.

\subsection{Excisive homology theories}\label{excthe}  
We consider triangulated categories  in terms of a loop functor $\Omega$.
A {\sf triangulated category}  $(\cT,\Omega, \mathcal{Q})$ is an additive
 category $\cT$ with an
 equivalence $\Omega:\cT \r \cT$ and a class $\mathcal{Q}$ of 
sequences  in $\cT$  called {\sf distinguished triangles}
$$
(T) \qquad \Omega C \r A \r B \r C 
$$
satisfying some axioms, see  \cite{Kel}, \cite{CT}. 
A {\sf triangle functor} from $(\cT_{1},\Omega_{1},\mathcal{Q}_{1})$ to
$(\cT_{2},\Omega_{2},\mathcal{Q}_{2})$ is a pair consisting of an additive functor
$R: \cT_{1}\r \cT_{2}$ and a natural transformation $\alpha:\Omega_{2}R \r R\Omega_{1}$
such that
$$
\Omega_{2}R(C)\xrightarrow{R(f)\circ \alpha_{C}}R(A) \xrightarrow{R(g)} R(B)\xrightarrow{ R(h)} R(C)
$$
is a distinguished triangle in $\cT_{2}$ for each triangle 
 $\Omega_{1}C\xrightarrow{f} A \xrightarrow{g} B \xrightarrow{h} C$.

Let $(\cT,\Omega, \mathcal{Q})$ be a triangulated category. 
An {\sf excisive homology theory}  for $\cC$ with values in $\cT$ is an $\mathcal{E}$-excisive homology, as defined in \cite[Sec 6.6]{CT},  when $\mathcal{E}$ is the class of weakly split extensions. 

Let $P$ be a property for functors defined on
$\cC$ to some triangulated category. 
A functor $u: \cC \r \cV$ is a {\sf
universal functor} with $P$ if it has the property $P$ and if $F:
\cC \r \cT$ is another functor with $P$ there exists a unique
triangle functor $G: \cV \r \cT$ such that the following
diagram is commutative
$$
\xymatrix{
\mbox{$\cC$} \ar[r]^-{u}\ar[dr]_{F}& \cV\ar@![d]^{G}\\
& \cT
}
$$
\subsection{The category $\mathfrak{KK}_{\cC}$} \label{catkk}  Let $A$ and $B$ be objects in $\cC$.
Consider $\cM_{\infty}$ the ind-ring defined in \cite[Sec 4.1]{CT} and define inductively
$B^{\cS^{n+1}}:=(B^{\cS^{n}})^{\cS^{1}}$.   Define
$$
E_{n}(A,B)_{\cC}:=[J^{n}(A),\cM_{\infty} B^{\cS^{n}}]_{\cC}.
$$
Consider the following morphism 
$\imath_{n}:E_{n}(A,B)_{\cC}\r E_{n+1}(A,B)_{\cC}$ such that
$$\begin{array}{c}
J^{n}(A)\xrightarrow{f}\cM_{\infty} B^{\cS^{n}}\\ 
\downarrow
\\
 J^{n+1}(A)
\xrightarrow{J(f)}
J(\cM_{\infty}B^{\cS^{n}})\xrightarrow{\phi_{\cM_{\infty},B^{\cS^{n}}}} 
\cM_{\infty}J(B^{\cS^{n}})\xrightarrow{\rho_{B^{\cS^{n}}}}\cM_{\infty}B^{\cS^{n+1}}
\end{array}
$$    
Define
$$
\kk_{\cC}(A,B)=\coli_{n\in \N}E_{n}(A,B)_{\cC}.
$$

Let $A$, $B$ and $C$ be objects in $\cC$. There exists an
associative product
$$
\circ:\kk_{\cC}(B,C)\times\kk_{\cC}(A,B)\r\kk_{\cC}(A,C)
$$
which extends the composition of algebra homomorphisms.
If $[\alpha]\in\kk_{\cC}(B,C)$ is an element represented by
$\alpha:J^{n}(B)\r C^{\cS^{n}}$ and $[\beta]\in \kk_{\cC}(A,B)$ is an
element represented by $\beta:J^{m}(A)\r B^{\cS^{m}}$ then
$[\alpha]\circ [\beta]$ is an element $\kk_{\cC}(A,C)$
represented by 
\begin{equation}\label{colaw}
J^{n+m}(A)\xrightarrow{J^{n}(\beta)}J^{n}(B^{\cS^{m}})\xrightarrow{}J^{n}(B)^{\cS^{m}}\xrightarrow{\alpha^{\cS^{m}}} C^{\cS^{n+m}}.
\end{equation}

This product allows us to define a composition in the category
$\mathfrak{KK}_{\cC}$ whose objects are the same objects of $\cC$ and the 
morphisms from $A$ to $B$ are the elements of $\kk_{\cC}(A,B)$.
Denote by 
$$
j_{\cC}:\cC \r \mathfrak{KK}_{\cC}
$$ 
the functor which at the level of objects is the identity and at
level of morphism sends $f:A\r B$ to $[f]\in \kk_{\cC}(A,B)$.
 \begin{rem}\label{redeid}
A morphism  $f:C\r M_{\infty} C$ in $\cC$
 represents an element $[f]$ in $\kk_{\cC}(C,\cM_{\infty}C)$. But
 also represents an element in $\kk_{\cC}(C,C)$ because  
$$
\kk_{\cC}(C,C)=\coli E_{n}(C,C)_{\cC} \quad\mbox{ and } \quad E_{0}(C,C)_{\cC}=[C,\cM_{\infty}C]_{\cC}.
$$ 
\end{rem}
 
Consider the functor
  $\Omega:\mathfrak{KK}_{\cC} \r \mathfrak{KK}_{\cC}$ 
 which sends an object 
  $A$ of $\cC$  to the path object $\Omega A$. Let $[\alpha]$ be an
  element of $\kk_{\cC} (A,B)$ represented by $\alpha: J^{n}(A)\r
  B^{\cS^{n}}$. The class of $\Omega [\alpha]$ is represented by
\begin{equation}\label{eq45}
J^{n}(A^{\cS^{1}})\xrightarrow\
J^{n}(A)^{\cS^{1}}\xrightarrow{\alpha^{\cS^{1}}} B^{\cS^{n+1}}
\end{equation}
Note \eqref{eq45} represents an element of
$\kk_{\cC}(A^{\cS^{1}},B^{\cS^{1}})$ (see Lemma 6.3.8 \cite{CT} to
check it is well defined). As $\imath:\Omega A \r A^{\cS^{1}}$ is a
$\kk_{\cC}$-equivalence (see corollary of Lemma 6.3.2 \cite{CT}), 
we have the natural isomorphism
$$ 
\kk_{\cC}(A^{\cS^{1}},B^{\cS^{1}})\xrightarrow{\sim}\kk_{\cC}(\Omega A, \Omega B).
$$
A diagram 
$$
\Omega C \r A\r B\r C
$$
of morphisms in $\mathfrak{KK}_{\cC}$ is called a {\sf distinguished
triangle} if it is isomorphic in $\mathfrak{KK}_{\cC}$ to the path sequence 
$$
\Omega B'\xrightarrow{j(\iota)} P_{f}\xrightarrow{j(\pi_{f})}
A'\xrightarrow{j(f)} B'
$$ 
associated with a homomorphism $f:A'\r B'$ in  $\cC$.
Denote the class of distinghished triangles by $\mathcal{Q}$.
The category $\mathfrak{KK}_{\cC}$ is triangulated with respect to
the  endofunctor 
$\Omega: \mathfrak{KK}_{\cC}\r\mathfrak{KK}_{\cC}$ and the class $\mathcal{Q}$ of
distinguished triangles, see \cite[Theorem 6.5.2]{CT}. 

Let $\operatorname{E}: A\xrightarrow{f}B\xrightarrow{g} C$ be a weakly split
extension and let 
$c_{\operatorname{E}}\in \kk_{\cC}(J(C),A)$ be the classifying map of $\operatorname{E}$.  
As the natural map $\rho_{A}:J(A)\r \Omega A$ induces a
$\kk_{\cC}$-equivalence (see Lemma 6.3.10, \cite{CT}) we can consider
the following morphisms in $\kk_{\cC}(\Omega C, A)$ 
\begin{equation}\label{ejekk}
\partial_{\operatorname{E}}:=c_{\operatorname{E}}\circ \rho^{-1}_{C}.
\end{equation}
The functor $j_{\cC}:\cC\r \mathfrak{KK}_{\cC}$ with the morphisms
 $\{\partial_{\operatorname{E}}:\operatorname{E}\in\cE\}$ is an excisive homology theory, 
homotopy invariant and $M_{\infty}$-stable.

\begin{thm}\label{unikk}
 The functor $j_{\cC}:\cC \r \mathfrak{KK}_{\cC}$ is universal
 with the properties defined above. In other words, if $\cT$
is a triangulated category and $G:\cC \r \cT$ together
a class of morphisms $\{\overline{\partial}_{\operatorname{E}}:\operatorname{E}\in\cE\}$ is an 
excisive, homotopy invariant and $M_{\infty}$-stable functor,
then there exists a unique triangle functor $\overline{G}:
 \mathfrak{KK}_{\cC}\r \cT$ such that the following diagram commutes
$$
\xymatrix{
\cC \ar[r]^{j_{\cC}}\ar[dr]_{G}& \mathfrak{KK}_{\cC}\ar@![d]^{\overline{G}}\\
&\cT
}
$$ 
\end{thm}
\proof 
By definition $j_{\cC}$ is homotopy invariant and $M_{\infty}$-stable. Let us show that $j_{\cC}:\cC \r \mathfrak{KK}_{\cC}$ is an excisive homology theory with the family  $\{\partial_{\operatorname{E}}:\operatorname{E}\in\cE\}$. Let $\operatorname{E}: A\xrightarrow{f}B\xrightarrow{g} C$ be a weakly split extension. Take the
 path sequence asociated to $g$ and the following diagram in $\mathfrak{KK}_{\cC}$
$$
\xymatrix{
\Omega C \ar[r]^{j(\iota)} & P_{g}\ar[r]^{j(\pi_{g})}& B \ar[r]^{j(g)}& C & T\\
\Omega C \ar[u]^{\id}\ar[r]_{\partial_{\operatorname{E}}} & A\ar[u]^{\iota_{f}}\ar[r]_{j(f)}&B\ar[u]_{\id}\ar[r]_{j(g)}&C\ar[u]_{\id}&T'
}
$$
The first square commutes beacuse $\iota_{f}\circ c_{\operatorname{E}}$ is elementarily
homotopic to $\iota\circ \rho_{C}$. By \cite[Lemma 6.3.2]{CT}, the morphism $\iota_{f}$ is a
$\kk_{\cC}$-equivalence. Finally $T$ and $T'$ are isomorphic in
$\mathfrak{KK}_{\cC}$ and $T'$ is a distinguished triangle. For the rest of the proof see \cite[Theorem 6.6.2]{CT}. \qed

Algebraic $K$-theory is not excisive, nor is it $M_{\infty}$-stable or homotopy invariant. However the homotopy algebraic $K$-theory, defined in \cite{kh} and denoted by $\kh$, have all these properties. Consider the funtor $\kh=\kh_{0}:\alg_{\ell}\r \operatorname{Ab}$. Because it satisfies \cite[Theorem 6.6.6]{CT} we have a natural map
$$
\kk_{\alg_{\ell}}(\ell,A)\r \hom_{\operatorname{Ab}}(\kh(\ell),\kh(A))
$$
As $\ell$ is unital, there is a map $\Z\r\ell$ which induce a map
$$
\hom_{\operatorname{Ab}}(\kh(\ell),\kh(A)) \mapsto \hom_{\operatorname{Ab}}(\kh(\Z),\kh(A))\simeq \kh(A). 
$$
Composing both maps we obtain a homomorphism $\kk_{\alg_{\ell}}(\ell, A)\r \kh(A)$ and the main theorem from \cite{CT} prove that it is an isomorphism. 
\begin{thm} \label{mainbak} 
Consider $\cC=\alg_{\ell}$ the category of $\ell$-algebras, then
$$
\kk_{\alg_{\ell}}(\ell, A)\simeq \kh(A).
$$
\end{thm}
\proof See \cite[Theorem 8.2.1]{CT}.\qed
\section{Equivariant matrix invariance}\label{eqmi}
In the equivariant setting we replace the property of
$M_{\infty}$-stability by a stability condition depending on $\cC$. 
We consider the different cases of $\cC$ separately.
\subsection{$G$-equivariant stability}
Regard $$M_{G}=\{f: G\times G \rightarrow \ell: \operatorname{supp}(f)<\infty \}$$ as the algebra of matrices with coefficients in $\ell$ indexed by $G\times G$,  with translation action of $G$: 
$$
g\cdot e_{s,t} = e_{gs,gt}.
$$

We are going to identify a $G$-algebra $A$ with the $G$-algebra $M_{G}\otimes A$ with the diagonal action. Note that the map $a\r a \otimes e_{1_{G},1_{G}}$ is not equivariant then we can not define $G$-stability as in Section \ref{maes}. 
For this reason we define $G$-stability in terms of $G$-modules.

A pair $(\mathcal{W},B)$ is a {\sf $G$-module with basis} if $\mathcal{W}$ is a
$G$-module, free as an $\ell$-module and $B$ is a basis of $\mathcal{W}$. A pair $(\mathcal{W}',B')$ is a {\sf submodule with basis} of
 $(\mathcal{W},B)$ if  $\mathcal{W}'$ is a submodule of $\mathcal{W}$ and $B'\subset B$. 
Note that if $(\mathcal{W}_{1},B_{1})$ and $(\mathcal{W}_{2},B_{2})$ are $G$-modules with basis
then $(\mathcal{W}_{1}\oplus \mathcal{W}_{2},B_{1}\sqcup B_{2})$ is a
$G$-module with basis. 

Let $(\mathcal{W},B)$ be a $G$-module with basis $B$. We define 
$$
\cL(\mathcal{W},B):=\{ \psi: B \times B\r \ell : \{v : \psi(v,w)\neq 0\}
\mbox{ is finite for all $w$}\}
$$ 
Note that  $\cL(\mathcal{W},B)$ and
$\en_{\ell}(\mathcal{W})=\{f:\mathcal{W}\r \mathcal{W} : \mbox{ $f$ is
$\ell$-linear}\}$ are isomorphic;
indeed we have inverse isomorphisms
\begin{equation}\label{bime}
\en_{\ell}(\mathcal{W}) \r \cL(\mathcal{W},B) \qquad f\mapsto \psi_{f} \qquad
\psi_{f}(v,w)=p_{v}(f(w))
\end{equation}
 here $p_{v}: \mathcal{W} \r \ell$ is the projection to the submodule of
 $\mathcal{W}$ generated by $v$
\begin{equation}\label{bime2}
 \cL(\mathcal{W},B) \r \en_{\ell}(\mathcal{W}) \qquad \psi \mapsto f_{\psi} \qquad
 f_{\psi}(w)=\sum_{v\in B} \psi(v,w) v.
\end{equation}
Define
$$
\begin{array}{l}
\cC(\mathcal{W},B):=\{\psi \in\cL(\mathcal{W},B):  \{w : \psi(v,w)\neq 0\}
\mbox{ is finite for all $v$}\}\\
\cF(\mathcal{W},B):=\{\psi \in\cL(\mathcal{W},B):  \{(v,w) : \psi(v,w)\neq 0\}
\mbox{ is finite }\}
\\
\en^{F}_{\ell}(\mathcal{W},B):=\{f\in \en_{\ell}(\mathcal{W}): \psi_{f}\in \cF(\mathcal{W},B)\}\\
\en^{C}_{\ell}(\mathcal{W},B):=\{f\in \en_{\ell}(\mathcal{W}): \psi_{f}\in \cC(\mathcal{W},B)\}
\end{array}
$$
Note that $\cC(\mathcal{W},B)$ is a ring with the matrix product and
$\en^{C}_{\ell}(\mathcal{W},B)$ is a ring with the composition. 
These rings are isomorphic.

Let $(\mathcal{W},B)$ be a $G$-module with basis. Consider the representation
$$
\rho: G\r \en_{\ell}(\mathcal{W}) \qquad\qquad \qquad \rho_{g}(w)=g\cdot w
$$
We say that $(\mathcal{W},B)$ is a {\sf $G$-module by finite automorphisms} if
 $\rho(G)\subset \en^{F}_{\ell}(\mathcal{W},B)$.
We say that $(\mathcal{W},B)$ is a {\sf $G$-module by locally finite automorphisms} if
 $\rho(G)\subset \en^{C}_{\ell}(\mathcal{W},B)$.
If $(\mathcal{W},B)$ is a $G$-module by locally finite automorphisms,
 $\en^{C}_{\ell}(\mathcal{W},B)$ and $\en^{F}_{\ell}(\mathcal{W},B)$ are $G$-algebras with the following action
$$
g\cdot f = \rho(g) f (\rho(g))^{-1}
$$
Note that $\en^{F}_{\ell}(\mathcal{W},B)$ is an ideal of $\en^{C}_{\ell}(\mathcal{W},B)$.

\begin{ex}\label{kgar}
Let $\mathcal{W}=\ell G$ be the group algebra considered as a
$G$-module via the regular representation with basis $B=\{\delta_{g}: g\in G\}$,
\begin{equation}
g\cdot (\sum_{h\in G}a_{h}\delta_{h})=\sum_{h \in G}a_{h}\delta_{gh}, \qquad a_{h}\in\ell.
\end{equation}
Note 
$$
\rho: G \r \en_{\ell}(\ell G)\simeq \cL(\ell G,B) \quad g\mapsto
M_{g}=\sum_{t\in G}e_{gt,t}  
$$
As $M_{g}\in \cC(\ell G,B)$ for all $g\in G$, $(\ell G,B)$ is a 
$G$-module by locally finite automorphisms. 
Moreover we have $(M_{g})^{-1}=M_{g^{-1}}=(M_{g})^{t}$ and 
$M_{G} = \cF(\ell G,B).$ 
\end{ex}

Let $A$ be a $G$-algebra. Consider the tensor product $M_{G} A =M_{G}\otimes A $
with the diagonal action of $G$.

\begin{rem}\label{196}
Let $\mathcal{W}$ be a
 $G$-module. Let $\mathcal{W}^{\tau}$ be $\mathcal{W}$ considered as
 a $G$-module with trivial action. Recall that $\ell G\otimes \mathcal{W} \simeq  \ell G \otimes \mathcal{W}^{\tau}.$
Inverse isomorphisms are given by
$$
\begin{array}{lll}
T:\ell G \otimes \mathcal{W}^{\tau} \r  \ell G\otimes \mathcal{W} &\quad\quad& S:\ell G \otimes \mathcal{W}
\r \ell G\otimes  \mathcal{W}^{\tau}\\
\\
T(\delta_{g}\otimes h) =\delta_{g}\otimes
g(h) &\quad\quad&
S(\delta_{g}\otimes h)=\delta_{g}\otimes g^{-1}(h)
\end{array}
$$
\end{rem}
\goodbreak
If $(\mathcal{W},B)$ is a $G$-module with basis we will write 
$\en^{C}_{\ell}(\mathcal{W})$ and $\en^{F}_{\ell}(\mathcal{W})$
ommiting the basis when there is no confusion.  
\begin{rem}\label{197}
Let $(\mathcal{W},B)$ be a $G$-module with basis. Let us check that
$$
\en^{F}_{\ell}(\ell G \otimes \mathcal{W} )\simeq \en^{F}_{\ell}(\ell G)\otimes\en^{F}_{\ell}(\mathcal{W}).
$$
Define a $G$-algebra homomorphism 
$T:\en^{F}_{\ell}(\ell G)\otimes\en^{F}_{\ell}(\mathcal{W})\r\en^{F}_{\ell}(\ell G\otimes
\mathcal{W})$ by
$$
T(e_{g,h}\otimes e_{v,w})=e_{g,v,h,w}\quad\quad v,w\in B\quad g,h\in G
$$
As $T$ is a bijection between the basis, $T$ is an isomorphism.
\end{rem}
Let $(\mathcal{W}_{1},B_{1})$ and $(\mathcal{W}_{2},B_{2})$ be $G$-modules by 
locally
finite auto\-morphisms such that 
$\operatorname{card}(B_{i})\leq\operatorname{card}(G)\times
 \operatorname{card}{\N}$, $i=1,2.$ 
The inclusion 
$\iota:\mathcal{W}_{1}\r \mathcal{W}_{1}\oplus \mathcal{W}_{2}$
induces a morphism of $G$-algebras 
\begin{equation}\label{itil}
\tilde{\iota} : \en^{F}_{\ell}(\mathcal{W}_{1})\r \en^{F}_{\ell}(\mathcal{W}_{1}\oplus \mathcal{W}_{2}) \qquad f \mapsto
\left(\begin{array}{ll}f & 0 \\ 0 & 0 \end{array}\right)
\end{equation}
Let $A$ be a $G$-algebra and consider 
$$
\tilde{\iota}\otimes 1: \en^{F}_{\ell}(\mathcal{W}_{1})\otimes A \r
 \en^{F}_{\ell}(\mathcal{W}_{1}\oplus \mathcal{W}_{2})\otimes A.
$$  
A functor $F:\mbox{$G$-$\alg$}\r \cD$ is {\sf $G$-stable} if  for
$(\mathcal{W}_{1},B_{1})$, $(\mathcal{W}_{2},B_{2})$ and $A$ as above
$F(\tilde{\iota}\otimes 1)$ is an isomorphism in $\cD$.

Let us show that if $F: \mbox{$G$-$\alg$}\r \Cat$ is a $G$-stable functor then $F$ is $M_{\infty}$-stable.
Consider $(\mathcal{W}_{1},B_{1})=(\ell,\{1\})$ and
$(\mathcal{W}_{2},B_{2})=(\ell^{(\N)}, \{e_{i}: i\in \N\})$ with 
$\ell^{(\N)}=\oplus_{i=1}^{\infty}\ell$, $\{e_{i}: i\in \N\}$ 
is the canonical basis and both modules have the trivial action of $G$.  
Note
$$
\en_{\ell}^{F}(\ell)=\en_{\ell}(\ell)=\ell \quad \en^{F}_{\ell}(\ell\oplus
\ell^{(\N)})=\en^{F}_{\ell}(\ell^{(\N)})=M_{\infty} 
$$
and $\tilde{\iota}_{\infty}: \ell \r M_{\infty}$ is the inclusion at the upper
left corner.
Then $\tilde{\iota}\otimes 1=\iota: A\r M_{\infty}(A)$ and $F(\iota)$
is an isomorphism. Observe that if $G=\{ e\}$, $F$ is $G$-stable if and only if $F$ is $M_{\infty}$-stable. 

Let $A$, $B$ be $G$-algebras and $F:\mbox{${G}$-$\alg$}\r \cD$ a functor. A {\sf zig-zag} between
$A$ and $B$ by $F$ is a diagram in ${G}$-$\alg$ 
$$
A\xrightarrow{f_{1}} C_{1} \xleftarrow{g_{1}}\ldots \xrightarrow{f_{n}} C_{n} \xleftarrow{g_{n}} B
$$ 
such that $F(g_{i})$, $i=1,\ldots, n$, are isomorphisms in $\cD$.

\begin{ex}\label{zzr}
Let $A$ be a $G$-algebra and $F$ a $G$-stable functor. There
exists a zig-zag between  $A$ and $M_{G}A$ by $F$. 
Consider $\mathcal{W}_{1}=(\ell G, B)$ as in the example \ref{kgar} and
consider $\mathcal{W}_{2}=(\ell, \{1\})$ with the trivial action of $G$.
Put $\mathcal{W}=\mathcal{W}_{1}\oplus \mathcal{W}_{2}$ and
$C=\en_{\ell}^{F}(\mathcal{W})$ 
with the induced action,
then
$$
\iota: A = A\otimes \ell= A\otimes \en_{\ell}^{F}(\ell) \r A\otimes C
\leftarrow A \otimes M_{G}: \iota'  
$$
is a zig-zag between $A$ and $M_{G}A$ by $F$.
\end{ex}
\begin{prop}\label{mesges}
Suppose $G$ is countable. Let $F:\mbox{$G$-$\alg$}\r \Cat$ be an $M_{\infty}$-stable functor. Define
$$
\hat{F}:\mbox{$G$-$\alg$}\r \Cat \quad \quad A\mapsto
F(M_{G}\otimes A)
$$  
Then $\hat{F}$ is $G$-stable.
\end{prop}
\proof
Let $(\mathcal{W}_{1},B_{1})$, $(\mathcal{W}_{2},B_{2})$ be $G$-modules by locally
finite automorphisms and let $A$ be a $G$-algebra. Consider 
$$
\tilde{\iota}\otimes 1: \en^{F}_{\ell}(\mathcal{W}_{1})\otimes A \r 
\en^{F}_{\ell}(\mathcal{W}_{1} \oplus \mathcal{W}_{2})\otimes A.
$$  
We have to prove that 
\begin{equation}\label{ecua}
\hat{F}(\tilde{\iota}\otimes 1):F(M_{G}\otimes\en^{F}_{\ell}(\mathcal{W}_{1})\otimes A) \r 
F(M_{G}\otimes \en^{F}_{\ell}(\mathcal{W}_{1} \oplus \mathcal{W}_{2})\otimes A)
\end{equation}
is an isomorphism. By remarks \ref{197} and \ref{196} we know that
$$
F(M_{G}\otimes \en_{\ell}^{F}(\mathcal{W}_{1})\otimes A)\simeq F( \en^{F}_{\ell}(\mathcal{W}_{1}^{\tau})\otimes M_{G} \otimes A)
$$
and
$$
F( M_{G}\otimes \en_{\ell}^{F}(\mathcal{W}_{1}\oplus \mathcal{W}_{2})\otimes A) \simeq F(\en^{F}_{\ell}((\mathcal{W}_{1}\oplus \mathcal{W}_{2})^{\tau})\otimes M_{G}\otimes A).
$$
Note that $\en_{\ell}^{F}(\mathcal{W}_{1}^{\tau})$ and $\en^{F}_{\ell}((\mathcal{W}_{1}\oplus \mathcal{W}_{2})^{\tau})$ are equivariantly isomorphic to $M_n$ or $M_\infty$.
As $F$ is $M_{\infty}$-stable, \eqref{ecua} is an isomorphism. \qed

\begin{rem}\label{ejpgj}
Suppose $G$ is a finite group of order $n$ such that $1/n\in \ell$. 
Then the element $\xi= (1/n)\sum_{g\in  G}\delta_{g}$ in $\ell G$ is idempotent. 
The map $s: \ell \r \ell G$, $s(1)=\xi$, is a $G$-equivariant section
of the canonical augmentation $\varphi: \ell G \r \ell$. Thus the
sequence of $G$-modules 
\begin{equation}\label{secex}
\xymatrix{
0\ar[r] & I \ar[r] & \ell G \ar[r]^{\varphi} &\ell \ar[r]& 0
} 
\end{equation}
splits.
Hence
$\ell G=\ell \xi \oplus I$. 
Notice that $I$ is a $G$-module with basis $\{\delta_{e}-\delta_{g}: g\neq e \}$. 
Define
$$
\lambda_{g}=\left\{\begin{array}{ll} \xi & g=e\\
\delta_{e}-\delta_{g}& g \neq e
\end{array}\right.
$$
The set  $\Lambda =\{ \lambda_{g}: g\in G\}$ is a basis of
$\ell G$ and the relations with the elements of
$B=\{\delta_{g}:g\in G\}$ are the following
$$
\begin{array}{lll}
\lambda_{e}=\frac{1}{n}\sum_{g\in G}\delta_{g} &\quad\quad&
\lambda_{h}=\delta_{e}-\delta_{h}\\ 
\delta_{e}=\lambda_{e} +
\frac{1}{n}\sum_{g\neq e}\lambda_{g} &\quad\quad&
\delta_{h}=\lambda_{e}-\lambda_{h}+\frac{1}{n}\sum_{g\neq e}
\lambda_{g}\end{array} \quad \quad h\neq e
$$
Consider $\mathcal{W}_{1}=\ell=(\ell\xi,
\{\xi\})$ and $\mathcal{W}_{2} = (I,\{\lambda_{g}\}_{g\neq e})$, in this case
the morphism 
$\eqref{itil}$ is
$$
\iota : \ell \r M_{G}\simeq \en({\ell G, \Lambda}) \quad 1\mapsto
\left(\begin{array}{llll}
1 &0 &\ldots &0\\
0&0&\ldots &0\\
\vdots& \vdots & \ddots&\vdots\\
0&0&\ldots &0
\end{array}\right) 
$$  
If we write it in the canonical basis $B$ we have
$$
\overline{\iota} : \ell \r M_{G}=\en({\ell G, B}) \quad 1\mapsto
\left(\begin{array}{llll}
\frac{1}{n} & \frac{1}{n}&\ldots &\frac{1}{n}\\
\frac{1}{n}&\frac{1}{n}&\ldots &\frac{1}{n}\\
\vdots& \vdots & \ddots&\vdots\\
\frac{1}{n}&\frac{1}{n}&\ldots &\frac{1}{n}
\end{array}\right) 
$$  
If $F:\mbox{${G}$-$\alg$}\r \Cat$ is a  $G$-stable functor then
$F(\overline{\iota})$ is an isomorphism in $\Cat$. 
\end{rem}

\subsection{$G$-graded stability}
In this section we consider a dual notion of $G$-equivariant stability. We want to identify a
$G$-graded algebra $A$ with the $G$-graded matrix algebra
 $M_{G}A$. The definition 
of $G$-graded stability is easier than that of $G$-equivariant stability  because the morphism $A\r M_{G}A$, $a\mapsto e_{1_{G}.1_{G}}\otimes a$, is homogeneous. 

Write $G_{gr}$-$\alg$ for the category of $G$-graded algebras.
Let $A$ be a $G$-graded algebra. Define the
following grading in $M_{G} A$ 
\begin{equation}{\label{gramg}}
(M_{G} A)_{g}:=<e_{s,t}\otimes a : g=s|a|t^{-1} >.
\end{equation}

Depending on the context $M_{G}A$  will be considered as an algebra, a $G$-algebra or a $G$-graded algebra.  
Let be
$$
\iota_{A}: A\r M_{G}A \qquad a \mapsto e_{1_{G},1_{G}}\otimes a,
$$ 
note $\iota_{A}$ is homogeneous.  A functor $F:\mbox{$G_{gr}$-$\alg$}\r \Cat$ is  {\sf $G_{gr}$-stable} if $F(\iota_{A})$ is an isomorphism in $\Cat$ for all $A\in \mbox{$G_{gr}$-$\alg$}$.

\begin{prop}\label{mesgesco}
Suppose $G$ is countable. Let $F:\mbox{$G_{gr}$-$\alg$}\r \Cat$ be an $M_{\infty}$-stable functor. Define
$$
\hat{F}:\mbox{$G_{gr}$-$\alg$}\r \Cat \quad \quad A\mapsto
F(M_{G} A)
$$  
with $M_{G} A$ as in \eqref{gramg}.
Then $\hat{F}$ is $G_{gr}$-stable.
\end{prop}
\proof
Denote by $M_{|G|}$ to $M_{G}$ with the trivial grading. Because $G$ is countable we have that $M_{\infty} M_{|G|}$ and $M_{\infty}$ are isomorphic in $G_{gr}$-$\alg$. We also have the following isomorphisms of $G$-graded algebras 
$$
\begin{array}{ccc}
\eta: M_{G}M_{G}\r M_{G}M_{|G|} &\quad &e_{s,g}\otimes e_{t,r} \mapsto e_{st,gr}\otimes e_{t,r}\\
\mu: M_{G}M_{|G|}\r M_{|G|}M_{G} &\quad & e_{s,g}\otimes e_{t,r} \mapsto e_{t,r}\otimes e_{s,g}
\end{array}
$$
Consider the following commutative diagram 
$$
\xymatrix{
& M_{G}M_{G}A\ar[r]^{\eta} &M_{G}M_{|G|}A \ar[dr]^{\mu}&\\
M_{G}A\ar[ur]^{\id_{M_{G}}\otimes \iota_{A}}\ar[dr]_{\iota_{\infty}}& & & M_{|G|}M_{G}A \ar[dl]^{\iota_{\infty}}\\
&M_{\infty}M_{G}A \ar[r]_-{\simeq}& M_{\infty}M_{|G|}M_{G}A&
}
$$
It follows that $\hat{F}(\iota_{A})=F(\id_{M_{G}}\otimes \iota_{A})$ is an isomorphism in $\Cat$ because $F$ is $M_{\infty}$-stable.
\qed

\section{Equivariant algebraic kk-theory}\label{ekk}
From now to the rest of the paper we suppose $G$ is a countable group.
\subsection{The category $\mathfrak{KK}^{G}$}\label{kkgcat}
In this section we introduce an equivariant algebraic kk-theory for $G$-algebras. 
 Let $A$, $B$ be $G$-algebras, we define 
$$
\kk^{G}(A,B):= \kk_{\mbox{\tiny{$G$-$\alg$}}} (M_{G}\otimes A, M_{G}\otimes B).
$$
Consider the category $\mathfrak{KK}^{G}$ whose objects are the
$G$-algebras and where the morphisms between $A$ and $B$ are the elements
of $\kk^{G}(A,B)$.  
Let  
$j^{G}:\mbox{$G$-$\alg$}\r \mathfrak{KK}^{G}$
be the functor defined as the identity on objects and which sends each
morphism of $G$-algebras $f:A\r B$ to its class $[\id_{M_{G}}\otimes f]\in
\kk^{G}(A,B)$.
The composition law in $\mathfrak{KK}^{G}$ is the same that in \eqref{colaw} taking $M_{G}\otimes A$, $M_{G}\otimes B$ and $M_{G}\otimes C$ instead of $A$, $B$ and $C$.
As in Section \ref{catkk} we have an equivalence  $\Omega: \mathfrak{KK}^{G}\r \mathfrak{KK}^{G}$ and distinguished triangles
$$
\Omega C \r A\r B\r C
$$
in $\mathfrak{KK}^{G}$ which gives to $\mathfrak{KK}^{G}$ a triangulated category structure. 

\begin{thm}\label{unikkg}
The functor $j^{G}:\mbox{$G$-$\alg$}\r \mathfrak{KK}^{G}$ is an
excisive, equivariantly homotopy invariant, and $G$-stable
functor. Moreover, it is the universal functor for these properties. In other words,  if $\cT$
is a triangulated category and $R:\mbox{$G$-$\alg$} \r \cT$ together
a class of morphisms $\{\overline{\partial}_{\operatorname{E}}:\operatorname{E}\in\cE\}$ is an 
excisive, equivariantly homotopy invariant and $G$-stable functor,
then there exists a unique triangle functor $\overline{R}:\mathfrak{KK}^{G}\r \cT$ such that the following diagram commutes
$$
\xymatrix{
 \mbox{$G$-$\alg$} \ar[r]^-{j^{G}}\ar[dr]_{R}& \mathfrak{KK}^{G}\ar@![d]^{\overline{R}}\\
&\cT
}
$$ 
\end{thm}
\proof Let $\operatorname{E}$ be a weakly split extension. Define 
$$
\partial^{G}_{\operatorname{E}}\in \hom_{\mathfrak{KK}^{G}}(\Omega C,
A)=\hom_{\mathfrak{KK}_{\mbox{\tiny{$G$-$\alg$}}}}(\Omega M_{G}\otimes C, M_{G}\otimes A)
$$
as the morphism $\partial_{\operatorname{E}'}$ defined in \eqref{ejekk} asociated to
the following weakly split extension
$$
M_{G}\otimes A\r M_{G}\otimes B\r M_{G}\otimes C \quad \quad (\operatorname{E}')
$$
By theorem \ref{unikk} and proposition \ref{mesges} the functor
$j^{G}:\mbox{$G$-$\alg$}\r \mathfrak{KK}^{G}$ with the family
$\{\partial^{G}_{\operatorname{E}}: \operatorname{E}\in\cE \}$ is an excisive, homotopy
invariant and $G$-stable functor.
Let us check it is universal for these properties. Let
$X:\mbox{$G$-$\alg$}\r \cT$ be a functor
which has the mentioned properties with a family 
$\{\overline{\partial}_{\operatorname{E}}:\operatorname{E}\in\cE \}$ . 
By theorem  \ref{unikk} there exists a
unique triangle functor $\overline{X}:\mathfrak{KK}_{\mbox{\tiny{$G$-$\alg$}}}\r\cT$
such that the following diagram commutes
\begin{equation}\label{warm}
\xymatrix{
\mbox{$G$-$\alg$}\ar[ddr]_-{X}\ar[rr]^{j^{G}}\ar[dr]^-{j_{\mbox{\tiny{$G$-$\alg$}}}}&& \mathfrak{KK}^{G}\ar@![ddl]^-{X'}\\
&\mathfrak{KK}_{\mbox{\tiny{$G$-$\alg$}}}\ar[ur]\ar[d]_-{\overline{X}}&\\
&\cT&\\
}
\end{equation}
We will define $X': \mathfrak{KK}^{G}\r\cT$. We know that
$X'=X$ on objects. 
As $X$ is $G$-stable the following morphisms are a zig-zag between
$A$ and $ M_{G}\otimes A$ by $X$,
(see Example \ref{zzr})
\begin{equation}
\xymatrix{
A \ar[r]^-{\iota_{A}}& M_{G\sqcup\{*\}}\otimes A & M_{G}\otimes A \ar[l]_-{\iota'_{A}}
}
\end{equation}
Let $\alpha \in \kk^{G}(A,B)$ and define 
$$
X'(\alpha):=X(\iota_{B})^{-1}X(\iota'_{B})\overline{X}(\alpha)X(\iota'_{A})^{-1}X(\iota_{A}).
$$ 
Note this definition is the unique possibility to make the diagram
\eqref{warm} commutative. 
\qed
\subsection{The category $\hat{\mathfrak{KK}}^{G}$}
In this section we introduce the equivariant algebraic kk-theory for $G$-graded algebras. 
Let $A$, $B$ be $G$-graded algebras. We define
$$
\hat{\kk}^{{G}}(A,B):= \kk_{\mbox{\tiny{$G_{gr}$-$\alg$}}}( M_{G} A, M_{G} B). 
$$
Consider the category $\hat{\mathfrak{KK}}^{G}$ whose objects are the
$G$-graded algebras and the morphisms between $A$ and $B$ are the
elements of $\hat{\kk}^{G}(A,B)$. 
Let $j^{G_{gr}}:\mbox{$G_{gr}$-$\alg$}\r \hat{\mathfrak{KK}}^{G}$
be the functor defined as the identity on objects and which sends each
morphism of $G$-graded algebras $f:A\r B$ to its class $[\id_{M_{G}} f]\in
\hat{\kk}^{G}(A,B)$. As in Section \ref{kkgcat} we can consider a triangulated category structure on $\hat{\mathfrak{KK}}^{G}$.

\begin{thm}\label{unikkhg}
The functor $j^{G_{gr}}:\mbox{$G_{gr}$-$\alg$}\r \hat{\mathfrak{KK}}^{G}$ is an
excisive, graded homotopy invariant, and $G$-graded stable functor. 
Moreover, it is the universal functor for these properties. In other words,  if $\cT$
is a triangulated category and $R:\mbox{$G_{gr}$-$\alg$} \r \cT$ together
a class of morphisms $\{\overline{\partial}_{\operatorname{E}}:\operatorname{E}\in\cE\}$ is an 
excisive, graded homotopy invariant and $G$-graded stable functor,
then there exists a unique triangle functor $\overline{R}:\hat{\mathfrak{KK}}^{G}\r \cT$ such that the following diagram commutes
$$
\xymatrix{
 \mbox{$G_{gr}$-$\alg$} \ar[r]^-{j^{G_{gr}}}\ar[dr]_{R}& \hat{\mathfrak{KK}}^{G}\ar@![d]^{\overline{R}}\\
&\cT
}
$$
\end{thm}
\proof The proof is similar to Theorem \ref{unikkg}.\qed

\section{Algebraic Green-Julg theorem}\label{tricrosec}
\subsection{Crossed product and trivial action}
Let $A$ be an $\ell$-algebra. Write $A^{\tau}$ for $A$ with the trivial action of  
$G$. This gives us a functor $\tau: \alg \r \mbox{$G$-$\alg$}$.
It is easy to check that $j^{G}\circ \tau$ satisfies excision, is homotopy
invariant and is $M_{\infty}$-stable. Write $\mathfrak{KK} = \mathfrak{KK}_{\alg}$.
 By Theorem \ref{unikk}
there exists a unique functor $\tau:\mathfrak{KK} \r
\mathfrak{KK}^{G}$ such that the following diagram is commutative
$$
\xymatrix{
\alg \ar[r]^-{\tau}\ar[d]_{j} & \mbox{$G$-$\alg$}\ar[d]^{j^{G}}\\
\mathfrak{KK}\ar[r]_-{\tau} &  \mathfrak{KK}^{G}
}
$$

Let $A$ be a $G$-algebra. The crossed product algebra $A \rtimes G$
is the $\ell$-module $A\otimes \ell G$ with the following multiplication 
$$
(a\rtimes g)(b\rtimes {h})=a(g\cdot b)\rtimes
{gh}\qquad a,b\in A\ g,h \in G.
$$
\begin{prop}\label{prcrsa}
Let $A$ be a $G$-algebra and $\mathcal{W}$ a $G$-module by locally finite
automorphisms. The following algebras are naturally isomorphic
$$
(A \rtimes G)\otimes \en^{F}_{\ell}(\mathcal{W})\simeq (A\otimes
\en^{F}_{\ell}(\mathcal{W}))\rtimes G.
$$  
\end{prop}
\proof
Let $\rho:G\r (\en^{K}_{\ell}(\mathcal{W}))^{\times}$ be the structure map.
Note that the homomorphisms 
$$
\phi:(A\rtimes G)\otimes \en^{F}_{\ell}(\mathcal{W})\r (A\otimes
\en^{F}_{\ell}(\mathcal{W})) \rtimes G\ \quad \quad 
\phi(a\rtimes{g}\otimes \varphi)=a\otimes
\varphi \rho(g^{-1})\rtimes {g}
$$
$$
\psi:(A\otimes
\en^{F}_{\ell}(\mathcal{W}))\rtimes G \r (A \rtimes G)\otimes
\en^{F}_{\ell}(\mathcal{W})\quad \quad 
\psi(a\otimes \varphi \rtimes {g})=a \rtimes
{g}\otimes\varphi \rho(g)
$$
are inverse of each other.   
\qed

\begin{prop}\label{crokk}
There exists a unique functor ${\rtimes} {G}:\mathfrak{KK}^{G}\r
\mathfrak{KK}$ such that the following diagram is commutative
$$
\xymatrix{
 \mbox{$G$-$\alg$}\ar[r]^-{\rtimes {G}}\ar[d]_{j^{G}} &\alg \ar[d]^{j}\\
\mathfrak{KK}^{G}\ar[r]_-{{\rtimes}{G}} &  \mathfrak{KK}
}
$$
\end{prop}
\proof
We shall show $j(- \rtimes {G})$ is excisive, homotopy invariant
and $G$-stable. Because
$\rtimes {G}$ maps split exact sequences to
split exact sequences and $j$ is excisive, then  $j(-
\rtimes {G})$ is excisive. 

That $j(-\rtimes G)$ is homotopy invariant follows from the fact that 
$$
A[t]\rtimes G = (A\rtimes G)[t].
$$
Let $(\mathcal{W}_{1},B_{1})$, $(\mathcal{W}_{2},B_{2})$ be
$G$-modules by locally finite 
automorphisms and $A$ a $G$-algebra. Consider the isomorphism $\psi$ 
defined in Proposition \ref{prcrsa}. 
Note that the following diagram is commutative
$$
\xymatrix{
 (A \otimes \en^{F}_{\ell}(\mathcal{W}_{1})) \rtimes G \ar[rr]^-{
(1\otimes \tilde{\iota})\rtimes G}\ar[d]_{\psi} && 
( A\otimes \en^{F}_{\ell}(\mathcal{W}_{1}\oplus \mathcal{W}_{2}))\rtimes G \ar[d]^{\psi}\\
(A \rtimes G)\otimes \en^{F}_{\ell}(\mathcal{W}_{1})\ar[rr]_-{(1\rtimes G)\otimes
\tilde{\iota}} && (A \rtimes G)\otimes \en^{F}_{\ell}(\mathcal{W}_{1}\oplus \mathcal{W}_{2}) 
}
$$
Because $j$ is $M_{\infty}$-stable, $j((1\rtimes G)\otimes
\tilde{\imath})$ is an isomorphism. Hence
$j(-\rtimes G)(1\otimes\tilde \imath)$ is an isomorphism by the
diagram above.\qed

\begin{rem}\label{reppc}
Let $[\iota_{\infty} \alpha]\in \kk^{G}(A,B)$ be an element represented by $\alpha:
J^{n}(M_{G}A)\r (M_{G}B)^{\operatorname{sd}^{p}\cS^{n}}$ which is a morphism in
 $[J^{n}(M_{G}A),\cM_{\infty}(M_{G}B)^{\operatorname{sd}^{p}\cS^{n}}]$.
Consider the classifying map
$$
J^{n}(M_{G}A\rtimes G)\r J^{n}(M_{G}A)\rtimes G
$$ 
The element $[\alpha]\rtimes {G}$ is represented by the following
composition
$$
J^{n}(M_{G}A\rtimes G)\r J^{n}(M_{G}A)\rtimes G
\xrightarrow{\alpha\rtimes G} (M_{G}B)^{\operatorname{sd}^{p}\cS^{n}}\rtimes G.
$$ 
\end{rem}
\begin{prop}
The functor ${\rtimes} {G}:\mathfrak{KK}^{G}\r \mathfrak{KK}$ is a triangle functor. 
\end{prop}
\proof A distinguished triangle in $\mathfrak{KK}^{G}$ is a diagram isomorphic to 
$$
\Omega B\xrightarrow{j^{G}(\iota)} P_{f}\xrightarrow{j^{G}(\pi_{f})}
A\xrightarrow{j^{G}(f)} B
$$
for some morphism of $G$-algebras $f:A\r B$. That means, if it is isomophic to 
\begin{equation}\label{trikg}
M_{G}\Omega B\xrightarrow{j_{\mbox{\tiny $G$-$\alg$}}(\iota)} M_{G}P_{f}\xrightarrow{j_{\tiny \mbox{$G$-$\alg$}}(\pi_{f})}
M_{G}A\xrightarrow{j_{\mbox{\tiny $G$-$\alg$}}(f)} M_{G} B
\end{equation}
in $\mathfrak{KK}_{\tiny \mbox{$G$-$\alg$}}$.  The functor ${\rtimes} {G}:\mathfrak{KK}^{G}\r \mathfrak{KK}$ sends the triangle \eqref{trikg} to 
\begin{equation}
(M_{G}\Omega B)\rtimes G\xrightarrow{} (M_{G}P_{f})\rtimes G \xrightarrow{} (M_{G}A)\rtimes G\xrightarrow{} (M_{G} B)\rtimes G
\end{equation}
which by Proposition \ref{prcrsa} and $M_{\infty}$-stability is isomorphic to 
\begin{equation}\label{estr}
\Omega (B\rtimes G)\xrightarrow{} P_{f}\rtimes G \xrightarrow{} A \rtimes G\xrightarrow{} B\rtimes G
\end{equation}
As $P_{f}\rtimes G\simeq P_{f\rtimes G}$, \eqref{estr} is a distinghished triangle in  $\mathfrak{KK}$.
\qed 
\subsection{Green-Julg Theorem for $\mathfrak{KK}^{G}$  }
In this section we shall see an algebraic version of
the {\sc Green-Julg Theorem}, see \cite{Mkk} and \cite{ght}  for versions of this result in Kasparov KK-theory and $E$-theory setting.   
\begin{thm}\label{gjgalg}
Let $G$ be a finite group of n elements and $1/n\in \ell$. 
The functors $\tau:\mathfrak{KK}
\r\mathfrak{KK}^{G}$  and
$\rtimes {G}:\mathfrak{KK}^{G} \r \mathfrak{KK}$ are adjoint functors.
Hence
$$
{\kk}^{G}(A^{\tau},B)\simeq \kk(A, B\rtimes G) \quad\quad A\in
\alg \quad \quad B\in \mbox{$G$-$\alg$}.
$$ 
\end{thm}
\proof
By \cite[Theorem 2, pag 81]{Mac}, it is enough to prove that there
exist natural transformations  
$\overline{\alpha}_{A}\in \kk(A, A^{\tau} \rtimes G)$ and $\overline{\beta}_{B}\in
\hkk((B \rtimes G)^{\tau},B)$ such that the following compositions 
 $$
 A^{\tau}\xrightarrow{\tau(\overline{\alpha}_{A})}
(A^{\tau}\rtimes G)^{\tau}\xrightarrow{\overline{\beta}_{\tau(A)}} A^{\tau} \qquad
 B \rtimes G \xrightarrow{\overline{\alpha}_{B\rtimes G}}
(B \rtimes G)^{\tau}\rtimes G \xrightarrow{\overline{\beta}_{B}\rtimes {G}} B\rtimes G
 $$
are the identities in $\hkk(A^{\tau},A^{\tau})$ and
$\kk(B\rtimes G, B \rtimes G)$ respectively.

Put $\epsilon=1/n\sum_{g\in G} {g}$ in
$\ell G$ and define 
\begin{equation}{\label{algj}}
\alpha_{A}:A\r A^{\tau}\rtimes G= A\otimes \ell G \quad\quad \alpha(a)=a\otimes \epsilon. 
\end{equation}
Note $\alpha_{A}$ is an algebra morphism since $\epsilon$ is idempotent.
Consider the element $\overline{\alpha}_{A} \in \kk(A, A^{\tau}\rtimes G)$ 
represented by $\alpha_{A}$.
Let
$$
\beta_{B}:(B\rtimes G)^{\tau}\r M_{G} B \quad\quad \beta_{B}(b\rtimes
 {g})=\sum_{s\in G}s(b)e_{s,sg}
$$
One can check that $\beta_{B}$ is an equivariant algebra morphism.
Let $\overline{\beta}_{B}\in \hkk((B\rtimes G)^{\tau},B)$
be represented by $\beta_{B}$.
The composite $\beta_{\tau(A)}\tau(\alpha_{A})$ is
$\id_{A^{\tau}}\otimes \overline{\iota}$ where 
$\overline{\iota}$ is the map defined in
\ref{ejpgj}. As $j^{G}$ is $G$-stable, $j^{G}(\id_{A^{\tau}}\otimes
\overline{\iota})$ is the identity in $\hkk(A^{\tau},A^{\tau})$.

Let $\psi$ be the morphism defined in the proof of the
Proposition  \ref{prcrsa}. By Remark \ref{reppc},
$\overline{\beta}_{B}\rtimes {G}$ is
represented by $\psi\circ (\beta_{B}\rtimes G)$.
We want to prove that $\psi\circ (\beta_{B}\rtimes G)\circ \alpha_{B\rtimes G}$ represents the identity in $\kk(B\rtimes G,
B\rtimes G)$. Note that 
$$
(\psi\circ (\beta_{B}\rtimes G)\circ \alpha_{B\rtimes G})(b\rtimes  {g}) = 
\frac{1}{n}\sum_{h,s\in G}(s(b)\rtimes  {h}) e_{s,h^{-1}sg}.
$$
Put $t=h^{-1}sg$ and note 
\begin{equation}\label{ecbbt}
\frac{1}{n}\sum_{h,s\in G}(s(b)\rtimes  {h}) e_{s,h^{-1}sg}= 
\frac{1}{n}\sum_{t,s\in G}(s(b)\rtimes  {sgt^{-1}}) e_{s,t}
\end{equation}
and
$$
s(b)\rtimes  {sgt^{-1}}=({\bf{1}}\rtimes {s})(b\rtimes
 {g})({\bf 1}\rtimes  {t^{-1}}) \quad\quad \mbox{ in
$\tilde{B}\rtimes G$}
$$
We can write \eqref{ecbbt} as $TA_{b\rtimes {g}}T^{-1}$, where
$$
A_{b\rtimes {g}}=\frac{1}{n}\sum_{t,s\in
  G}(b\rtimes {g})e_{s,t} \quad T=\sum_{t\in G}({\bf 1}\rtimes  {t})e_{t,t} 
$$
Because $b\rtimes {g}\mapsto
A_{b\rtimes {g}}$ represents the identity, the same is true for  
$b\rtimes {g}\mapsto TA_{b\rtimes {g}}T^{-1}$, see \cite[Proposition 5.1.2]{CT}. \qed

\begin{ex} 
We give an example to show that the adjointness between of $\tau$ and
$\rtimes {G}$ of Theorem \ref{gjgalg} fails to hold at the algebra level.
Let $G=\Z_{2}=\{1,\sigma\}$, $A=\ell$ and $B=(\ell G)^{*}$ the dual
algebra of $\ell G$ with the regular action.   
Note $\hom_{\mbox{$G$-$\alg$}}(A^{\tau},B)$  has two elements only:
$$
\varphi_{i}:\ell \r (\ell G)^{*} \qquad \varphi_{0}(1)=0 \qquad \varphi_{1}(1)=\chi_{1}+\chi_{\sigma}
$$
One the other hand $\hom_{\alg}(A, B\rtimes
G)=\hom_{\alg}(\ell,(\ell G)^{*}\rtimes G)$ 
has at least as many elements as $\ell$.
For each $\lambda \in \ell$ we can define 
$$
\varphi_{\lambda}:\ell \r(\ell G)^{*}\rtimes G \quad \varphi_{\lambda}(1)=
\chi_{1}\rtimes  {1}+\lambda(\chi_{1}\rtimes  {\sigma}) \quad
\lambda \in \ell
$$
Note $\varphi_{\lambda}$ is an algebra morphism because
$\chi_{1}\rtimes  {1}+\lambda(\chi_{1}\rtimes  {\sigma})$ is
an idempotent element.

\end{ex}
Write  $\psi_{GJ}$ for the isomorphism of the Theorem 
\ref{gjgalg}  
\begin{equation}\label{pgj}
\psi_{GJ}: \hkk(B^{\tau},A)\r \kk(B,A \rtimes G) \quad\quad
\psi_{GJ}=\alpha^{*}\circ \rtimes {G} 
\end{equation}
where $\alpha$ is the morphism defined in \eqref{algj}.

\begin{cor}\label{corogj} 
Let $G$ be a finite group such that $1/|G| \in \ell$. Let $A$ be a
$G$-algebra, then 
$$
\hkk(\ell,A)\simeq \kk(\ell, A \rtimes G)\simeq \kh(A \rtimes G)
$$
\end{cor}\qed

\section{Induction and Restriction}\label{inressec}
 \numberwithin{equation}{section}

In this section we study the adjoitness property of the functors
of induction and restriction  between 
$\mathfrak{KK}^{H}$ and $\mathfrak{KK}^{G}$ where $G$ is a group and $H$ is a subgroup of $G$. 

Let $A$ be a $G$-algebra and $H\subset G$ a subgroup. 
If we restrict the action to $H$ we obtain
an $H$-algebra $\resg(A)$. It is clear this
construction defines a functor $\resg: G$-$\alg \r
H$-$\alg$. 
It is easily seen that we can extend $\resg: G$-$\alg \r
H$-$\alg$ to a triangle functor ${\resg}: \mathfrak{KK}^{G}\r
\mathfrak{KK}^{H}$ so that the following diagram commutes
$$
\xymatrix{
\mbox{$G$-$\alg$}\ar[rr]^-{\resg}\ar[d]_{j^{G}}&&
\mbox{$H$-$\alg$}\ar[d]^{j^{H}}\\
\mathfrak{KK}^{G}\ar[rr]_-{{\resg}}&& \mathfrak{KK}^{H}
}.
$$

Let $\pi:G\r G/H$ the projection and
$A$ an $H$-algebra. Consider
$$
A^{(G,H)}:=\{f: G \r A : \#\pi(\operatorname{supp}(f))<\infty \}
$$
and define 
$$
\indg(A)=\{ f\in A^{(G,H)}: \quad f(s)=h(f(sh)) \quad \forall s\in
G, \quad h\in H \}.
$$
One checks that $\indg(A)$ is a $G$-algebra with pointwise
multiplication and the following action of $G$
\begin{equation}\label{acind}
(g\cdot f)(s)=f(g^{-1}s) \quad f\in \indg(A) \quad g,s\in G.
\end{equation}
Observe this construction is functorial, if $\varphi: A\r B$  is a morphism of $H$-algebras then 
$\indg(\varphi)(f)=\varphi \circ f$.

If $g\in G$, write $\chi_{g}:G \r \Z$ for the characteristic
function. If $a\in A$ and $g\in G$, define
\begin{equation}\label{defxi}
\xi_{H}(g,a)=\sum_{h\in H} \chi_{gh}h^{-1}(a) \qquad
\xi_{H}(g,a)(s)=\left\{
\begin{array}{ll}
h^{-1}(a) & s=gh\\
0 & s \notin  gH 
\end{array}
\right.
\end{equation}
It is easy to check that these elements belong to $\indg(A)$ and every element $\phi\in \indg(A)$ can be written as a finite sum
\begin{equation}\label{uniquexp}
\phi=\sum_{g\in \mathcal{R}}\xi_{H}(g,\phi(g))
\end{equation}
where $r:G/H\r G$ is a pointed section and $\mathcal{R}=r(G/H)$. Note that we have the following relations
\begin{equation}\label{actixi}
s\cdot \xi_{H}(g,a)=\xi_{H}(sg,a)
\end{equation}
\begin{equation}\label{prodxi}
\xi_{H}(g,a)\xi_{H}(\tilde{g},\tilde{a})= \left\{ \begin{array}{ll}
\xi_{H}(\tilde{g},\tilde{g}^{-1}g(a)\tilde{a}) & \tilde{g}^{-1}g\in H\\
0&  \tilde{g}^{-1}g\notin H
\end{array}\right. 
\end{equation}
\begin{equation}\label{hxi}
\xi_{H}(g,a)=\xi_{H}(gh,h^{-1}\cdot a)\quad h\in H
\end{equation}

\begin{prop}\label{proirc}
Let $A$ be a $G$-algebra and $B$ be an $H$-algebra, then
$$
\indg(B\otimes \resg A)\simeq \indg(B)\otimes A
$$
\end{prop}
\proof
The isomorphisms are given by 
$$
\begin{array}{c}
\begin{array}{llll}
S:&\indg(B)\otimes A &\r &\indg(B\otimes \resg (A))\\
&\xi_{H}(g, b)\otimes a &\mapsto & \xi_{H}(g,b\otimes g^{-1}\cdot a)
\end{array}
\\
\\
\begin{array}{llll}
T:&\indg(B\otimes \resg A) &\r &\indg(B)\otimes A\\
&\xi_{H}(g,b\otimes a) &\mapsto & \xi_{H}(g,b)\otimes g\cdot a
\end{array}
\end{array}
$$
\qed
\begin{cor}
Let $A$ be a $G$-algebra. Then
$$
\indg\resg A \r \ell^{(G/H)}\otimes A \qquad \qquad
\xi_{H}(s,b)\mapsto \chi_{sH}\otimes s\cdot b
$$
is an isomorphism of $G$-algebras.
\end{cor}\qed

\begin{prop}
Let $\operatorname{Ind}:H$-$\alg\r G$-$\alg$ be the following functor
$$\operatorname{Ind}(A)=\indg(M_{H}\otimes A).$$ There exists a functor
 ${\indg}:\mathfrak{KK}^{H}\r\mathfrak{KK}^{G}$ such that
 the following diagram is commutative
\begin{equation}\label{diaind}
\xymatrix{
\mbox{$H$-$\alg$}\ar[rr]^-{\operatorname{Ind}}\ar[d]_{j^{H}}&&
\mbox{$G$-$\alg$}\ar[d]^{j^{G}}\\
\mathfrak{KK}^{H}\ar[rr]_-{{\indg}}&& \mathfrak{KK}^{G}
}
\end{equation}
\end{prop}
\proof
Straightforward.\qed

\begin{prop}\label{indtri}
The functor ${\indg}:\mathfrak{KK}^{H}\r\mathfrak{KK}^{G}$ is a triangle functor.
\end{prop}
\proof Let $f: A\r B$ be a morphism of $H$-algebras. The following is an isomorphism of $G$-algebras 
\begin{equation}\label{yoquese}
\Theta: \indg(P_{f})\r P_{\indg(f)} \qquad \xi_{H}(g,(tp(t),a))\mapsto (\xi_{H}(g,tp(t)),\xi_{H}(g,a)).
\end{equation}
The image of the induction functor applied to the path extension of $f$ is
\begin{equation}\label{larala}
\indg (\Omega B)\xrightarrow{\indg \iota} \indg(P_{f})\xrightarrow{\indg \pi_{f}}\indg  A\xrightarrow{\indg f} \indg B
\end{equation}
By Proposition \ref{proirc} and \eqref{yoquese}, the extension \eqref{larala} is isomorphic to
$$
\Omega \indg ( B)\xrightarrow{\indg \iota} P_{\indg f} \xrightarrow{ \pi_{\indg f}}\indg  A\xrightarrow{\indg f} \indg B.
$$
Then ${\indg}:\mathfrak{KK}^{H}\r\mathfrak{KK}^{G}$ is a triangle functor.
\qed

\begin{thm}\label{adjinres}
Let $G$ be a  group and $H$ a subgroup of $G$. Then the functors
$$
\indg : \mathfrak{KK}^{H}\r  \mathfrak{KK}^{G} \quad \quad \resg :
 \mathfrak{KK}^{G}\r  \mathfrak{KK}^{H} 
$$
are adjoint. Hence
$$
\hkk(\indg(B),A)\simeq {\kk}^{H}(B,\resg(A))\quad \quad \forall
B\in \mbox{$H$-$\alg$} \quad A\in \mbox{$G$-$\alg$}
$$
\end{thm}
\proof
Let $A\in G$-$\alg$ and $B\in H$-$\alg$.  We need  natural
transformations 
$$
\alpha_{A}\in \hkk(\indg\resg A, A) \qquad
\beta_{B}\in {\kk}^{H}(B,\resg \indg B) \quad
$$ 
which verify the unit and counit condition, respectively.

Define $\varphi_{A}:\indg(\resg(A))\r M_{G/H}\otimes A$ such that
$$
\varphi_{A}(\xi_{H}(s,b))=e_{sH,sH}\otimes s\cdot b.
$$
One checks that $\varphi_{A}$ is a $G$-equivariant algebra morphism. 
Put
\begin{equation}\label{pside}
\psi_{B}: B \r \resg\indg(B)\quad \psi_{B}(b)=\xi_{H}(e,b)
\end{equation}
It is easy to check that $\psi_{B}$ is well-defined and is a map 
of $H$-algebras.
Let $\alpha_{A} \in \hkk(\indg\resg A, A)$  the element represented by
$\varphi_{A}$ and $\beta_{B}\in {\kk}^{H}(B,\resg \indg B)$ the
element represented by $\psi_{B}$. 

The composite $\resg (\alpha_{A})\circ \beta_{\resg A}$
is represented by $\resg(\varphi_{A})\circ\psi_{\resg A}$
which is ${\kk}^{H}$-equivalent to the identity in the sense of
remark \ref{redeid}.
The element $\alpha_{\indg B}\circ \indg (\beta_{B}) \in \hkk(\indg B, \indg B)$
is represented by 
$$
\gamma: \indg B\r  M_{G/H}\otimes \indg B\qquad \xi_{H}(g,b)\mapsto e_{gH,gH}\otimes \xi_{H}(g,b) 
$$
The following morphism of $H$-algebras 
$$
\theta: C\r M_{G/H}\otimes C \qquad\quad  \theta(c)=e_{H,H}\otimes c
$$
represents to the identity in the sense of remark \ref{redeid}.
Then $\indg(\theta)$ is $\hkk$-equivalent to the identity. It is easy
to check $\indg(\theta)=\gamma$ with $C= B$.
\qed 

Write $\psi_{IR}$ for the isomorphism 
\begin{equation}{\label{pir}}
\psi_{IR}: \hkk(\indg B,A) \r {\kk}^{H}(B,\resg A) \quad \quad
\psi_{IR}={\psi_{B}}^{*}\circ \resg
\end{equation}
where $\psi_{B}$ is the morphism defined \eqref{pside}.

\begin{cor}
Let $G$ be a group, $H$ a finite subgroup of $G$ and  $A$ a
$G$-algebra  then
$$
\hkk(\ell^{(G/H)},A)\simeq \kk(\ell, A\rtimes H ) \simeq
\kh(A\rtimes H)
$$
\end{cor}
\proof The isomorphism is the composition of 
$\psi_{GJ}$ and $\psi_{IR}$ defined in \eqref{pgj} and in \eqref{pir}.\qed

\section{Baaj-Skandalis Duality}\label{dualsec}

In this section we define crossed product functors between the categories $G$-$\alg$ and 
$G_{gr}$-$\alg$. We prove that they extend to equivalences 
between $\mathfrak{KK}^{G}$  and $\hat{\mathfrak{KK}}^{G}$.
In this way we obtain an algebraic duality theorem similar to the
duality  given by Baaj-Skandalis in \cite{bs}.

Let $A$ be a $G$-algebra. Then 
$$
A \rtimes G =\bigoplus_{s\in G} A\rtimes  {s} \qquad \mbox{
and } (A\rtimes
 {s})(A\rtimes  {t})\subset A\rtimes  {st}
$$
thus $A\rtimes G$ is a $G$-graded algebra.  
If $f: A\r B$ is a homomorphism of $G$-algebras then $f\rtimes
G: A\rtimes G \r B\rtimes G$ is a graded homomorphism. Hence we
have a  functor
$$
\rtimes {G}: G\mbox{-$\alg$}\r G_{gr}\mbox{-$\alg$}
$$
We can also define a functor
$$
{G}\hat{\rtimes}: G_{gr}\mbox{-$\alg$}\r G \mbox{-$\alg$}
$$
as follows. Let $B$ be a $G$-graded algebra. Let
 $G\hat{\rtimes} B$ be the algebra
which as a module is $\ell^{(G)}\otimes B$ and the product is the following
\begin{equation}\label{hgpro}
(\chi_{g}\rtimes a)(\chi_{h}\rtimes b):=\chi_{g}\rtimes a_{g^{-1}h}b.
\end{equation}
Here $b_{g}$ is the homogeneous element associated to
$g$ in the decomposition
$$
b=\sum_{g\in G} b_{g}.
$$ 
One checks that the product \eqref{hgpro} is associative and the
 action of $G$, $s\cdot \chi_{g}\rtimes a=\chi_{sg}\rtimes a$, makes it into a $G$-algebra.
The crossed product $G\hat{\rtimes} B \rtimes G$ is a $G$-graded algebra which contains $B$ as a graded subalgebra by
\begin{equation}
b\mapsto \sum_{g\in G}\chi_{1_{G}}\rtimes b_{g}\rtimes g.
\end{equation}
In \eqref{gramg} we defined $M_{G}B$ a $G$-graded matrix algebra asociated to $B$. Through the inclusion $b\mapsto e_{1_{G},1_{G}}\otimes b$ we can see $B$ as a $G$-graded subalgebra of it.  
In Proposition \ref{prole} below, we prove that $G\hat{\rtimes} B \rtimes G$ is isomorphic to $M_{G} B$ as a $G$-graded algebras. 

If $f:A\r B$ is a homogeneous homomorphism we define a $G$-algebra homomorphism $(G\hat{\rtimes} f)$ in the obvious way. 
Thus we have a functor
\begin{equation}{\label{funy}}
\hat{\rtimes} {G}: \mbox{$G_{gr}$-$\alg$}\r \mbox{$G$-$\alg$}.
\end{equation}

\begin{prop}\label{prole}
Let $A$ be a $G$-algebra and let $B$ be a $G$-graded algebra. 
\begin{itemize}
\item [a)] There are natural isomorphisms of $G$-algebras 
$$
G\hat{\rtimes} (A \rtimes G) \simeq M_{G}\otimes A
$$
\item [b)] There are natural isomorphisms of $G$-graded algebras
$$
(G \hat{\rtimes} B) \rtimes G \simeq M_{G}\otimes B
$$
\end{itemize}
\end{prop}
\proof
\begin{itemize}
\item [a)]
Define $T: G\hat{\rtimes} (A \rtimes G) \r M_{G}\otimes A$ as
$$T(\chi_{g}\rtimes a \rtimes  {s})= g\cdot a \otimes e_{g,gs}.$$
It is easy to check that $T$ is an equivariant algebra isomorphism with 
inverse given by
$$
S( a \otimes e_{r,t}):= \chi_{r}\rtimes r^{-1}\cdot a\rtimes
 {r^{-1}t} $$

\item[b)]
Define $T : (G\hat{\rtimes} B)\rtimes G \r M_{G}\otimes B$ as 
$$
T(\chi_{h}\rtimes b \rtimes
 {s})=\sum _{r\in G}e_{h,s^{-1}hr}\otimes b_{r}.
$$
It is easy to check that $T$ is a graded algebra isomorphism with inverse
given by
\begin{equation}\label{mpd}
S(e_{r,s}\otimes b_{q})= \chi_{r}\rtimes b_{q}\rtimes  {rqs^{-1}}
\end{equation}\qed
\end{itemize}

\begin{thm}\label{algbsd}\mbox{}
The functors $\rtimes G$ and $G\hat{\rtimes}$ extend to inverse
equivalences 
$$
-\rtimes G  :\mathfrak{KK}^{G} \longrightarrow
 \hat{\mathfrak{KK}}^{G} \quad \quad
 G \hat{\rtimes}-:\hat{\mathfrak{KK}}^{G} \longrightarrow 
 \mathfrak{KK}^{{G}}
$$
Hence if $A$ and $B$ are $G$-algebras and $C$ and $D$ are
$G$-graded algebras then
$$
\kk^{G}(A,B)\simeq \hat{\kk}^{G}(A\rtimes G, B\rtimes G)\quad 
\hat{\kk}^{G}(C,D)\simeq \kk^{G}( G\hat{\rtimes} C, G\hat{\rtimes}
D)$$
\end{thm}
\proof 
As $\rtimes {G}$ maps split sequences to
split sequences, $j^{G}(-\rtimes {G})$ is
excisive. By Proposition \ref{crokk}  $j^{G}(-\rtimes {G})$ 
is $G$-stable and homotopy invariant, whence it extends to 
${-\rtimes {G}}:\mathfrak{KK}^{G}\r
\hat{\mathfrak{KK}}^{G}$ by universality.
Similary, as $G\hat{\rtimes}$ maps split exact sequences to
split exact sequences then $j^{G}(G\hat{\rtimes}-)$ is
excisive. Because  $G\hat{\rtimes}$ maps graded homotopies to
equivariant homotopies and $j^{G}({G\hat{\rtimes} -})$ is
$M_{\infty}$-stable, $j^{G}(G\hat{\rtimes}-)$ extends to
$G \hat{\rtimes} -:\hat{\mathfrak{KK}}\r \mathfrak{KK}^{G}$
by universality.
To finish we must show that the maps 
$$
\kk^{G}(A,B)\r \hat{\kk}^{G}(A\rtimes G, B\rtimes G)
\quad \mbox{ and } \quad \hat{\kk}^{G}(C,D)\r
\kk^{G}(G\hat{\rtimes} C, G\hat{\rtimes} D) 
$$
are isomorphisms. This is true by Proposition \ref{prole}.
\qed

\bibliographystyle{plain}

\end{document}